\documentclass[10pt,oneside]{article}

\usepackage{graphicx}

\makeatother

\normalsize
\usepackage{latexsym}
\usepackage{amsfonts}       % ,epsfig}
\usepackage{amsmath}
\usepackage[mathscr]{euscript}
\usepackage{color}
\newtheorem{theo}{Theorem}[section]
\newtheorem{lemma}[theo]{Lemma}
\newtheorem{prop}[theo]{Proposition}
\newtheorem{cor}[theo]{Corollary}

\newcommand{\be}{\begin{equation} \label}       %
\newcommand{\ee}{\end{equation}}                %
\newcommand{\bea}{\begin{eqnarray}\label}
\newcommand{\eea}{\end{eqnarray}}                %
\newcommand{\bas}{\begin{eqnarray*}}            %
\newcommand{\eas}{\end{eqnarray*}}              %
\begin{document}
\title{Constructing solutions for a kinetic model of
angiogenesis in annular domains}
\author{ A. Carpio, G. Duro, M. Negreanu}
\date{January 23, 2016}
\maketitle
\begin{abstract}
We prove existence and stability of solutions for a model of angiogenesis 
set in an annular region. Branching and extension of blood vessel tips is described by
an integrodifferential kinetic equation of Fokker-Planck type supplemented with nonlocal boundary conditions and coupled to a diffusion problem with Neumann boundary conditions through the force field created by the tumor induced 
angiogenic factor and the flux of vessel tips.  
Our technique exploits balance equations, estimates of velocity decay and compactness results for kinetic operators, combined with gradient estimates of heat kernels for Neumann problems  in non convex domains.
\end{abstract}

\section{Introduction}
\label{sec:intro} 

Angiogenesis is a process through which new blood vessels grow from pre-existing ones. Angiogenesis is vital for tissue delevopment and repair. However, angiogenic disorders are often the cause of inflammatory and immune diseases \cite{angiogenesis1}. Moreover, angiogenesis is essential for the transition of benign tumors into malignant ones, and for subsequent tumor spread \cite{angiogenesis1}.
Numerous antitumor therapies target blood vessel growth \cite{angiogenesis2} in an attempt to prevent tumor expansion. Mathematical models may help to control the formation and evolution of blood vessel networks for therapeutical purposes. Many models have been proposed to describe different aspects of the process, see references \cite{angiogenesis3,angiogenesis6,angiogenesis5,angiogenesis4} for instance. However, the incessant availability of new experimental observations promotes  continued model update and fosters the search for improved descriptions. 

\begin{figure}
\centering
\includegraphics[width=7cm]{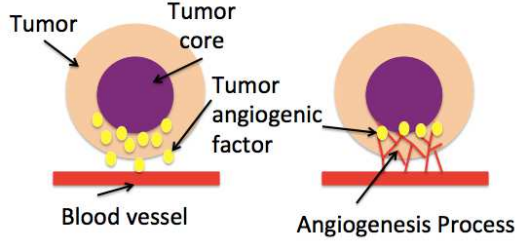}
\caption{Schematic representation of the formation of a vessel network to increase oxygen supply towards the inner regions of a tumor from a neighboring blood vessel. Drops represent the emitted concentration of tumor angiogenic factor, decreasing from the tumor core in the direction of the closest vessel tips.}
\label{fig1}
\end{figure}

In a tumor induced angiogenic process, high cell density in the
inner regions of the tumor results in low oxygen and nutrient levels. Cells respond emitting a substance (the tumor angiogenic factor) that eventually reaches neighboring blood vessels, promoting the appearance of new vessel tips that advance in direction to the tumor to supply new resources to the necrotic cells, see Fig. \ref{fig1}. 
The stochastic evolution of the vessel branching process seems to be a key feature to be taken into account. Recently, a deterministic integrodifferential system has been shown to  reproduce some aspects of the development of the stochastic vessel network \cite{capasso}.  
The evolution of the density of blood vessel tips $p$ in response to the concentration of tumor angiogenic factor released by cells $c$ is described by the following set of equations:
\begin{eqnarray} \frac{\partial}{\partial t} p(\mathbf{x},\mathbf{v},t)&=&
 \alpha(c(\mathbf{x},t)) \delta_{\mathbf v_0}(\mathbf{v}) p(\mathbf{x},\mathbf{v},t)  - \gamma p(\mathbf{x},\mathbf{v},t) \int_0^t d\,s \int d{\bf v}' p(\mathbf{x},\mathbf{v}',s)  \nonumber\\
 && - \mathbf{v}\cdot \nabla_\mathbf{x}   p(\mathbf{x},\mathbf{v},t) 
 + \beta {\rm div}_\mathbf{v} (\mathbf{v} p(\mathbf{x},\mathbf{v},t))+
\nonumber\\
&& - {\rm div}_\mathbf{v} \left[\mathbf{F}\left(c(\mathbf{x},t))\right)p(\mathbf{x},\mathbf{v},t)  \right]\! + \sigma
\Delta_\mathbf{v} p(\mathbf{x},\mathbf{v},t), \label{eq:p}  \\
\frac{\partial}{\partial t}c(\mathbf{x},t) &=& d \Delta_{\mathbf x} c(\mathbf{x},t) - \eta c(\mathbf{x},t) j(\mathbf{x},t) \label{eq:c}, \\
p(\mathbf{x},\mathbf{v},0) &=& p_0(\mathbf{x},\mathbf{v}), \quad
c(\mathbf{x},0) =c_0(\mathbf{x}), \label{eq:pc0}
\end{eqnarray}
where
\begin{eqnarray}
\alpha(c)=\alpha_1\frac{\frac{c}{c_R}}{1+\frac{c}{c_R}}, \quad
 {\bf F}(c)= \frac{d_1}{(1+\gamma_1c)^{q_1}}\nabla_{\mathbf x} c,
\label{eq:alphaF} \\
j(\mathbf{x},t)= \int_{\mathbb R^N} 
{|\mathbf{v}| \over 1 + e^{|\mathbf{v}- \mathbf{v}_0 \chi |^2 / \sigma_v^2}}
p(\mathbf{x},\mathbf{v},t)\, d \mathbf{v},
\quad \rho(\mathbf{x},t)= \int_{\mathbb R^N} p(\mathbf{x},\mathbf{v},t)
\, d \mathbf{v}, \label{eq:intpintvp}
\end{eqnarray}
for ${\mathbf x} \in \Omega \subset  \mathbb{R}^N$, ${\mathbf v} \in  \mathbb{R}^N$, 
$N=2, 3$, $t \in [0, \infty).$ The constants $\beta$, $\sigma$, $\gamma$, $d$, $\eta$, $\alpha_1$, $c_R$, $d_1$,  $\gamma_1$, $q_1$ are positive. The parameter $\chi >>1$ (typically $\chi >10$) whereas $\sigma_v^2 <<1.$ 
$\delta_{\mathbf v_{0}}$ is a Dirac measure supported at a point $\mathbf v_0$.
$\mathbf v_0$ is a typical sprouting velocity for the tips.  
The source term  $\alpha(c) \delta_{\mathbf v_0} \, p$ represents creation of new tips due to vessel tip branching. Tip vessel death when a tip encounters another vessel (anastomosis) is described by the integral sink  $- \gamma p \int_0^t \rho(p)$.  The Fokker-Planck operator expresses blood vessel extension.  
The chemotactic force $\mathbf F(c)$ is taken to depend on the flux of blood vessel tips through $j$ to represent that consumption of tumor angiogenic factor  is mostly due to the additional endothelial cells that produce vessel extensions \cite{angiogenesis3}. The velocity cut-off through the Fermi-Dirac distribution in the definition of $j$ (\ref{eq:intpintvp}) reflects the fact that cell velocities are limited, and small \cite{unbounded}.

We study the existence of solutions to a regularized version of equations 
(\ref{eq:p})-(\ref{eq:intpintvp}), where $\delta_{\mathbf v_0}$ is approximated 
by a smooth, positive, integrable and bounded function $\nu(\mathbf v)$, 
\begin{eqnarray} \frac{\partial}{\partial t} p(\mathbf{x},\mathbf{v},t)&=&
 \alpha(c(\mathbf{x},t)) \nu(\mathbf{v}) p(\mathbf{x},\mathbf{v},t)  - \gamma p(\mathbf{x},\mathbf{v},t) \int_0^t ds \, \rho(\mathbf{x},s)  \nonumber\\
 && - \mathbf{v}\cdot \nabla_\mathbf{x}   p(\mathbf{x},\mathbf{v},t) 
 + \beta {\rm div}_\mathbf{v} (\mathbf{v} p(\mathbf{x},\mathbf{v},t))+
\nonumber\\
&& - {\rm div}_\mathbf{v} \left[\mathbf{F}\left(c(\mathbf{x},t))\right)p(\mathbf{x},\mathbf{v},t)  \right]\! + \sigma
\Delta_\mathbf{v} p(\mathbf{x},\mathbf{v},t), \label{eq:pnu}
\end{eqnarray}
when $\Omega$  is an annular domain $r_0 < r < r_1$. 
Notice that delta functions can be approximated by sequences
of gaussians. The motivation for the annular geometry is simple, in view of Figure \ref{fig1}. Many tumors resemble spheres. An inner necrotic core is surrounded by a corona through which blood vessels spread, driven by the tumor angiogenic factor emitted by core. New vessel tips arise from existing vessels surrounding the outer layers of the tumor. They spread to supply with blood inner tumor regions in need of oxygen and nutrients. 

The general form of the boundary conditions in dimension $N=2,3,$ is as follows. We impose Neumann boundary conditions for $c$:
\begin{eqnarray}
{\partial c \over \partial r}({\mathbf x},t) = c_{r_0}({\mathbf x},t)<0, 
\quad {\mathbf x} \in S_{r_0},  \quad
 {\partial c \over \partial r} ({\mathbf x},t) = 0, \quad
{\mathbf x} \in S_{r_1},  \quad t \in [0,T],  \label{bc:c} 
\end{eqnarray}
where $c_{r_0}$ represents the influx of tumor angiogenic factor coming from
the inner core of the tumor.  $S_{r_0}$ and $S_{r_1}$ are spheres of
radius $r_0$ and $r_1$, respectively.

Since diffusion is absent in the $\mathbf x$ variable, the transport operator 
forces boundary conditions of the form:
\begin{equation}
p^-(\mathbf{x},\mathbf{v},t) = g(\mathbf{x},\mathbf{v},t) 
\quad {\rm on} \; \Sigma_T^-. \label{bc:plin}
\end{equation}
The sets $\Sigma_T^{\pm} =  (0,T)\times \Gamma^{\pm}, $ where
$
\Gamma^{\pm} = \{ (\mathbf{x},\mathbf{v}) \in \partial \Omega 
\times \mathbb{R}  \; |\; \pm \mathbf{v} \cdot \mathbf{n}(\mathbf{x})  >0\},
$
$\mathbf{n}(\mathbf{x})$ being the outward unit normal onto the boundary
$\partial \Omega$.  We denote by $p^+$ and $p^-$ the traces of $p$
on $\Sigma_T^{+}$ and $\Sigma_T^{-}$, respectively. In our geometry,
the boundary conditions for $p$ are defined using 
the magnitudes that can actually be measured: the marginal tip density 
$\rho = \int p d\mathbf v$ in the inner boundary and the flux of blood vessels
$\mathbf j = \int \mathbf v p d \mathbf v$ in the outer boundary. 
Using coordinates $\mathbf x= r \boldsymbol \theta$, with $r= |\mathbf x|,$
$\boldsymbol \theta \in S_{N-1}$, and
$\mathbf v= v_r \boldsymbol \phi$, with $v_r= |\mathbf v|$, $\boldsymbol \phi \in S_{N-1}$, the boundary conditions on $\Sigma_T^-$ read:
\begin{eqnarray}
 p^-(r_0, \boldsymbol \theta, v_r, \boldsymbol \phi,t) =  
 {e^{-{\beta \over \sigma} |\mathbf v \!-\! \mathbf v_0|^2} \over {\cal I}_0} 
  \Big[ \rho(r_0,\! \boldsymbol \theta, \!t) \!-\!
 \int_0^\infty  \hskip -4mm  d  \tilde v_r  \tilde v_r^{N-1}  \hskip -2mm
\int_{ \{\tilde {\boldsymbol \phi} \in S_{N\!-\!1} | \tilde{\mathbf v} 
\cdot \mathbf n  > 0 \} } 
 \hskip -22mm d \tilde   {\boldsymbol \phi} \, 
 p^+(r_0,  \! \boldsymbol \theta, \! \tilde v_r,  \! \tilde {\boldsymbol \phi}, \!t) 
 \Big],  \label{bc:Nr0} \\
 p^-(r_1,\boldsymbol \theta,v_r, \boldsymbol \phi,t) = 
 {e^{-{\beta \over \sigma} |\mathbf v \!-\! \mathbf v_0|^2} \over {\cal I}_1}
 \Big[ \!-\! j_0 \!-\!  \int_0^\infty \hskip -4mm
d \tilde v_r   \tilde v_r^{N-1} \hskip -2mm
\int_{ \{ \tilde {\boldsymbol \phi} \in S_{N\!-\!1} | 
\tilde {\mathbf v} \cdot \mathbf n  >0 \} } 
\hskip -22mm
d \tilde {\boldsymbol \phi} \, p^+(r_1, \! \boldsymbol \theta, \! \tilde v_r, 
\!\tilde {\boldsymbol \phi}, \!t) f_1(\mathbf v)  \Big], \label{bc:Nr1}
\end{eqnarray}
where $p^+$ and $p^-$ denote the traces of the solution $p$ on
$\Sigma_T^+$ and $\Sigma_T^-$, respectively, and
\begin{eqnarray}
{\cal I}_0= \int_0^\infty \hskip -4mm d \tilde v_r  \tilde v_r^{N-1}  
\hskip -2mm
\int_{ \{ \tilde {\boldsymbol \phi} \in S_{N\!-\!1} | 
\tilde {\mathbf v} \cdot \mathbf n  < 0 \} }
\hskip -20mm d \tilde {\boldsymbol \phi} \;
e^{-{\beta \over \sigma} |\tilde{\mathbf v} \!-\! \mathbf v_0|^2}, \quad
{\cal I}_1= \int_0^\infty  \hskip -4mm d \tilde v_r  \tilde v_r^{N-1}  \hskip -2mm
\int_{ \{ \tilde {\boldsymbol \phi} 
\in S_{N\!-\!1} | \tilde {\mathbf v} \cdot \mathbf n  < 0 \} }
\hskip -20mm
d \tilde \phi \; e^{-{\beta \over \sigma} |\tilde{\mathbf v} \!-\! \mathbf v_0|^2}
f_1(\tilde  {\mathbf v}).
\label{bc:integralsN}
\end{eqnarray}
The remaining functions are defined as:
\begin{eqnarray}
f_1(\mathbf v) = \mathbf v \cdot \mathbf n \Big[1 + 
e^{|\mathbf v - \mathbf v_0 \chi|^2  / \sigma_v^2} \Big]^{-1},
\label{bc:f1} \\
j_0(\boldsymbol \theta,t)=v_0\,\alpha(c(r_1,\boldsymbol \theta,t))\, 
p(r_1,\boldsymbol\theta,v_0,\mathbf w_0,t), 
\label{bc:j0N} 
\end{eqnarray}
for the fixed  velocity $\mathbf v_0 = (v_0,\mathbf w_0,)$, $v_0>0$, $\mathbf w_0 \in \mathbb R^{N-1}$.
Notice that the operators defining these boundary conditions are positive.
Thus, these conditions are expected to be absorbing, for positive densities. Similar boundary conditions are employed in kinetic models of charge transport in semiconductors \cite{holger}.

Rigorous derivations of these mean field models from the original stochastic systems as well as the development of stable numerical schemes require well posedness results for the integrodifferential set of equations. 
Equation (\ref{eq:p}) evokes Vlasov-Poisson-Fokker-Planck (VPFP) systems, with several key differences. First, the force field ${\bf F}$ is not related to the marginal tip density $\rho(p)$ through a Poisson equation. It depends on the flux of vessel tips ${j}$ through the gradient of solutions of heat equations with Neumann boundary conditions. Second, it contains a quadratic anastomosis term involving a 
nonlocal in time integrodifferential sink.  Moreover, the structure of the boundary conditions for the transport operator differs from those usually considered in Boltzmann equations for gas dynamics \cite{cercignani,hamdache} and
studied for VPFP models \cite{carrillo} as well, see also references \cite{chen, velazquez}.
Existence results for VPFP systems and related models in the whole space have been formulated under successively milder assumptions, see references
\cite{degond,victoryclassical,victoryweak,rein,bouchut,diperna,perthame}.
Global solutions for this angiogenesis model in the whole space have been constructed in  \cite{unbounded,unboundedheat}. 
Spatial boundaries pose new difficulties, arising from the nonlocal boundary conditions for the transport operator in the equation for the density of blood 
vessel tips and the presence of Neumann boundary conditions in the diffusion equation for the tumor angiogenic factor. 
Analyses in unbounded domains rely heavily on the properties of fundamental solutions for linear operators. The unavailability of results on fundamental solutions in bounded domains forces the development of new strategies.

In this paper, we prove existence and stability  of solutions of regularized versions of (\ref{eq:c})-(\ref{bc:plin}) where the measure $\delta_{\mathbf v_0}$ is replaced by a smooth positive bounded function.  Solutions are constructed as limits of solutions of linearized problems where all the nonlocal coefficients, rather than the sink terms, are frozen. This guarantees the nonnegativity of the densities $p$ and concentrations $c$, but requires $L^\infty_{\mathbf x}$ estimates of velocity integrals.  
Controlling the velocity decay of the densities provides such estimates.
Comparison principles and integral inequalities for both the diffusion and the kinetic equation allow us to control the $L^q$ norms of their solutions. Energy arguments provide basic derivative estimates. To handle the nonlocal coupling of the Neumann problem  with the kinetic equation we will have to make use of the theory of heat kernels in bounded domains \cite{heatkernel,varadhan1,varadhan2} and sharp gradient estimates for the semigroup of the Neumann problem \cite{chino} established by probabilistic methods for non convex regions in order to obtain $L^r-L^q$ estimates of the derivatives of solutions. Compactness results specific of kinetic operators \cite{perthame,diperna,bouchut}
enable the passage to the limit in the linearized problems.

The paper is organized as follows. In Section \ref{sec:linear} we adapt existence,
uniqueness and stability results for linear boundary value problems involving
Fokker-Planck operators, introducing additional lower order terms.
Section \ref{sec:decayv} derives $L^{\infty}$ estimates for the nonlocal coefficients
defined as velocity integrals of the vessel tip densities.  Bounds
on the velocity decay are essential to pass to the limit in linearized iterative
schemes that freeze the nonlocal coefficients.
In Section  \ref{sec:heat} we study the Neumann problem set in the annulus, 
establishing sharp estimates on the gradient of the solutions. These bounds are
fundamental to control the force field created by the tumor angiogenic factor. 
Section \ref{sec:nonlinearbc} proves the existence and stability result for the 
nonlinear problem with fixed known boundary condition. Finallly, Section
\ref{sec:nonlocalbc} addresses the angiogenesis problem with nonlocal
boundary conditions.

\section{Boundary value problems for linear Fokker-Planck operators}
\label{sec:linear}

Solutions for the coupled angiogenesis model will be constructed using an iterative
scheme that uncouples and freezes each variable to update the other.
A good knowledge about the properties of solutions of uncoupled linearized
equations is essential. In this section, we collect the needed existence
results and estimates for our specific linear problem for the density.

Let $\Omega \subset \mathbb{R}^N$ be a $C^\infty$ bounded domain 
with boundary $\partial \Omega$. Let us introduce the set
$Q_T= \Omega\times \mathbb{R}^N\times(0,T)$, $T>0$.
We consider the problem:
 \begin{eqnarray}
 \frac{\partial p}{\partial t} + \mathbf v \cdot\nabla_{\mathbf x} p
 +{\rm div}_{\mathbf v}(({\bf F}-\beta \mathbf v)p)-\sigma\Delta_{\mathbf v} p+ap=h 
  &\mbox{in } &  Q_T, 
\label{lin:ec}  \\
p({\mathbf x}, {\mathbf v},0)=p_0({\mathbf x}, {\mathbf v})  
  & \mbox{on } &   \Omega \times\mathbb{R}^N,  
\label{lin:in}
\end{eqnarray}
with $\beta\geq 0$, $\sigma>0$, 
${\bf F}({\mathbf x},t) \in L^\infty(\Omega\times(0,T))^N $ and $a \in L^\infty(Q_T)$. 
We will encounter two typical situations:
\begin{itemize}
\item $h\geq 0$, $a\in L^\infty(\Omega\times(0,T))$, $a \geq 0$, 
\item $h=0$, $a=a^+ - a^-$, $a^+ \in L^\infty(\Omega\times(0,T))$, $a^-\in 
L^\infty(Q_T)$.
\end{itemize}
The initial state $p_0$ represents a density. Therefore, $p_0 \geq 0.$ 
The transport operator selects absorbing boundary conditions of the form 
(\ref{bc:plin}) with $g \geq 0$.  

We seek weak solutions (in distributional sense) of the problem.  For any $T>0,$ a function $f\in L^\infty(0,T;L^1(\Omega\times \mathbb{R}^N))$ is a weak solution of equations (\ref{lin:ec})-(\ref{lin:in}) with boundary condition 
(\ref{bc:plin}) if  
\begin{equation} \label{lin:weak}
\hskip -2mm
 \begin{array}{cc}
 \displaystyle
\int_{Q_T} \hskip -3mm p \left[\frac{\partial \varphi}{\partial t}
+ \mathbf v\cdot \nabla_{\mathbf x} \varphi
-\beta {\mathbf v}\cdot \nabla_{\mathbf v} \varphi
+ {\bf F}\cdot \nabla_{\mathbf v} \varphi
+\sigma\Delta_{\mathbf v} \varphi-a\varphi\right] \,d{\mathbf x} \,d{\mathbf v}\, dt \\ 
\displaystyle
+\int_{\Omega\times \mathbb{R}^N} \hskip -6mm p_0
\varphi({\mathbf x},{\mathbf v},0) \,d{\mathbf x} \,d{\mathbf v} 
+\int_{\Sigma^-_T} \hskip -2mm 
|{\mathbf v}\cdot {\mathbf n}({\mathbf x})|g\varphi 
\, dS \, d{\mathbf v} \,dt
=\int_{Q_T} \hskip -2mm h\varphi \,d{\mathbf x} \,d{\mathbf v}\, dt
\end{array}  \end{equation}
for any $\varphi\in C^\infty_0(\overline{\Omega}\times \mathbb{R}^N\times[0,T))$ such that $\varphi=0$ on $\Sigma^+_T$.

We denote by $L^q$ the standard spaces of functions $p$ for which
$|p|^q$ is integrable with respect to the Lebesgue measure in
the pertinent domains and  by $L^{\infty}$ the space of bounded functions.
We introduce the space $L^q_k(\Sigma^{\pm}_T)$ of functions $g$ such that $|g|^q $ is integrable in $\Sigma^{\pm}_T$ with respect to the kinetic measure 
$|{\bf v} \cdot {\bf n}({\bf x})|  dS d{\mathbf v} dt$, where $dS$ is the Lebesgue measure on $\partial\Omega$. 
In an analogous way, we define $L^q_k(\Gamma^{\pm})$ with respect to the
measure $|{\bf v} \cdot {\bf n}({\bf x})|  dS d{\mathbf v}$. 

In absence of the lower order term $ap$, existence, smoothness, positivity and uniqueness results were established in reference \cite{carrillo}. Most of them 
extend  to the case $a\neq 0$  with slight modifications to the proofs. 

\begin{theo} {\bf (Existence, uniqueness, positivity).}  
\label{positivity}
Let $\Omega \subset \mathbb{R}^N$ be a bounded domain and set $T>0$. If
\begin{itemize}
\item[i)] ${\mathbf F} \in L^{\infty}(\Omega \times (0,T))$, 
$a \in L^{\infty}(Q_T),$
\item[ii)] $h \in L^2(Q_T)$, $p_0 \in L^2(\Omega \times   \mathbb{R}^N)$  and $g \in L^2_k(\Sigma^-_T),$ 
\end{itemize}
there exists a unique solution $p$ of equations (\ref{lin:ec})-(\ref{lin:in}),
(\ref{bc:plin}), satisfying:
\begin{itemize}
\item $p \in \{ f\in L^2(Q_T) \; | \;  \frac{\partial}{\partial t} f
+ \mathbf{v}\cdot \nabla_\mathbf{x}  f -\beta  \mathbf{v} \cdot \nabla_\mathbf{v} f
\in L^2(Q_T) \}.$
\item The equations hold in the sense of distributions: for any
$\phi \in C^{\infty}_0(\overline \Omega\times \mathbb{R}^N\times [0,T))$ and any $T>0$
\begin{eqnarray}
\int_{Q_T} p \left(
{\partial \phi \over \partial t} + {\mathbf v} \cdot \nabla_{\mathbf x} \phi
- \beta {\mathbf v} \cdot \nabla_{\mathbf v} \phi 
+ {\mathbf F} \cdot \nabla_{\mathbf v}  \phi
+ \sigma \Delta_{\mathbf v} \phi - a \phi
\right)  \,d{\mathbf x} \,d{\mathbf v} \,dt   \nonumber \\
+ \int_{\Omega \times  \mathbb{R}^N} \hskip -4mm
 p_0 \phi({\mathbf x},{\mathbf v},0) \,d{\mathbf x}\, d{\mathbf v}
 = \int_{\Sigma_T} \hskip -2mm
  {\mathbf v} \cdot {\mathbf n}({\mathbf x})  \, {\rm Tr} p\,
 \phi  \,dS \,d{\mathbf v} \,dt +
 \int_{Q_T} \hskip -2mm h\varphi \,d{\mathbf x} \,d{\mathbf v}\, dt.
  \label{weak}
\end{eqnarray}
If $\phi= 0$ on $\Sigma^+_T$, the boundary integral becomes
$$ -\int_{\Sigma^-_T} ({\mathbf v} \cdot {\mathbf n}({\mathbf x})) \, g\,
 \phi \,dt \,dS \,d{\mathbf v}.$$
\item ${\rm Tr} \,p=g$ on $\Sigma^-_T$ and $p(\mathbf{x},\mathbf{v},0) = 
p_0(\mathbf{x},\mathbf{v})$ in $ \Omega \times \mathbb{R}^N$.
\item 
$\|p\|_{L^{\infty}(0,T;L^2(\Omega \times  \mathbb{R}^N))} 
\leq C_1 \left[  \|p_0\|_{L^2(\Omega \times \mathbb{R}^N)} \!+\!
\|g\|_{L^2_k(\Sigma^-_T)} \!+\! \|h\|_{L^2(Q_T)}  \right],$
where $C_1>0$ depends on $T$, $\beta$ and $\|a^-\|_{\infty}$.
\item If $h\geq 0$, $p_0\geq 0$ and $g\geq 0$, then $p\geq 0$.
\end{itemize}
\end{theo}

\begin{theo} {\bf (Smoothness, balance laws).}  
\label{conservation}
Let $\Omega \subset \mathbb{R}^N$ be a bounded domain and set $T>0$. If
\begin{itemize}
\item[i)] ${\mathbf F} \in L^{\infty}(\Omega \times (0,T))$, $a \in L^{\infty}(Q_T),$
\item[ii)] $h \in L^1 \cap L^{\infty}(Q_T)$  and $|{\mathbf v}|^2 h \in L^{\infty}(0,T;L^1(\Omega \times \mathbb{R}^N))$,
\item[iii)] $p_0 \in L^1 \cap L^{\infty}(\Omega \times   \mathbb{R}^N)$ 
and $|{\mathbf v}|^2 p_0  \in L^1(\Omega \times    \mathbb{R}^N)$,
\item[iv)] $g \in L^1_k \cap L^{\infty}_k(\Sigma^-_T)$ and
$|{\mathbf v}|^2 g \in L^1_k(\Sigma^-_T)$,
\end{itemize}
the solution $p$ of equations (\ref{lin:ec})-(\ref{lin:in}),(\ref{bc:plin}) satisfies
\begin{itemize}
\item $p \in L^{\infty}(0,T;L^1\cap L^{\infty}(\Omega \times \mathbb{R}^N))$,
\item $|{\mathbf v}|^2 p \in L^{\infty}(0,T;L^1(\Omega \times \mathbb{R}^N))$,
\item $\nabla_{\mathbf v} p \in L^2(Q^T)$ and ${\rm Tr} \,p \big|_{\Sigma^+_T} \in
L^2_k(\Sigma^+_T)\cap L^{\infty}(0,T;L^1_k(\Gamma^+))$,
\item ${\rm Tr} \, p^2 \big|_{\Sigma^+_T} \in 
L^{\infty}(0,T;L^1_k(\Gamma^+))$,
\item Balance of mass: The solution $p$ has trace values in $L^\infty(0,T; L^1_k(\Gamma{+}))$ and verifies the continuity equation in integral form 
\begin{eqnarray}
{d\over dt} \int_ {\Omega \times \mathbb{R}^N}   \!\!\!\!\!\!\!\! p \,d{\mathbf x}\, d{\mathbf v}
&=& \!\! \int_{\Gamma^-} \!\!\!\! |{\mathbf v} \cdot {\mathbf n}({\mathbf x}) |
g \,dS\, d{\mathbf v} +\int_ {\Omega \times \mathbb{R}^N}   \!\!\!\!\!\!\!\! h \,d{\mathbf x}\, d{\mathbf v}  \label{mass} \\
&-& \!\! \int_{\Gamma^+} \!\!\!\! |{\mathbf v} \cdot {\mathbf n}({\mathbf x}) |
{\rm Tr} \,p \,dS\, d{\mathbf v} -  \int_ {\Omega \times \mathbb{R}^N} \!\!\!\!\!\!\!\! a p \,d{\mathbf x}\, d{\mathbf v}, \nonumber
\end{eqnarray}
\item Balance of momentum: If $|{\mathbf v}|^\mu h \in L^{\infty}(0,T;L^1(\Omega \times \mathbb{R}^N))$ and $|{\mathbf v}|^\mu g \in L^{\infty}(0,T;L^1_k(\Gamma^-))$, then $m_\mu(p)= \int_ {\Omega \times \mathbb{R}^N}   |{\mathbf v}|^\mu p \,d{\mathbf x}\, d{\mathbf v}$ is absolutely continuous and
\begin{eqnarray}
{d\over dt} \int_ {\Omega \times \mathbb{R}^N} \!\!\!\!\!\!\!\!\!\! 
|{\mathbf v}|^\mu  p \,d{\mathbf x}\, d{\mathbf v}  = 
 \!\int_{\Gamma^-} \!\!\! |{\mathbf v} \cdot {\mathbf n}({\mathbf x}) | |{\mathbf v}|^\mu g \,dS\, d{\mathbf v}
+\int_ {\Omega \times \mathbb{R}^N}   \!\!\!\!\!\!\!\! |{\mathbf v}|^\mu h \,d{\mathbf x}\, d{\mathbf v} \label{momentum}\\
 - \!\int_{\Gamma^+}\!\!\!\!\! |{\mathbf v} \!\cdot\!  {\mathbf n}({\mathbf x}) | 
 |{\mathbf v}|^\mu {\rm Tr} \,p \,dS\, d{\mathbf v} -\beta \mu \!\int_ {\Omega \times \mathbb{R}^N}\!\!\!\!\!\!\!\!\!\!  |{\mathbf v}|^\mu  p \,d{\mathbf x}\, d{\mathbf v}  
 -  \!\int_ {\Omega \times \mathbb{R}^N}\!\!\!\!\!\!\!\!\!\! 
a  |{\mathbf v}|^{\mu}  p \,d{\mathbf x}\, d{\mathbf v} \nonumber \\
+ \mu(\mu-2+N) \sigma  \!\!\int_ {\Omega \times \mathbb{R}^N}\!\!\!\!\!\!\!\!\!\! 
|{\mathbf v}|^{\mu-2}  p \,d{\mathbf x}\, d{\mathbf v}
 + \mu \!\int_ {\Omega \times \mathbb{R}^N}\!\!\!\!\!\!\!\!\!\!  
{\mathbf F} \cdot {\mathbf v} |{\mathbf v}|^{\mu-2}  p \,d{\mathbf x}\, d{\mathbf v},
\nonumber
\end{eqnarray}
\item $L^q$ estimates: If $h\geq 0$, $g\geq 0$ and $p_0\geq0$, then $p\geq 0$ 
and
\begin{eqnarray} \hskip -8mm
\frac{d}{d t}\|p(t)\|^q_{L^q(\Omega\times \mathbb{R}^N)} \hskip -3mm &=&
\hskip -3mm
\int_{\Gamma{-}}|{\mathbf v}\cdot {\mathbf n}({\mathbf x})|g^q dS \, d{\mathbf v}
+ q \int_{\Omega \times \mathbb R^2} 
\hskip -4mm h p^{q-1} d{\mathbf x} \, d{\mathbf v}    \nonumber \\
&-&  \hskip -3mm \int_{\Gamma^{+}}|\mathbf v\cdot \mathbf n(\mathbf x)|
{(\rm Tr} p)^q  dS \, d{\mathbf v}
-q \int_{\Omega \times \mathbb R^2} 
\hskip -4mm a p^q d{\mathbf x} \, d{\mathbf v}   
\nonumber \\
&+& \hskip -3mm N\beta(q\!-\!1)\|p(t)\|^q_{L^q} 
\hskip -2mm
- \sigma q(q\!-\!1) \!\!\int_{\Omega\times \mathbb{R}^N} \hskip -6mm 
p^{(q\!-\!2)}|\nabla_{\mathbf v} p|^2  d{\mathbf x} \, d{\mathbf v},
\label{lp}
\end{eqnarray}
% Correcto, sale todo dividido por q menos el gradiente y la a y h y al % multiplicar quedan las otras sin el q que divide y esa con el q
for any $1\leq q<\infty$. Setting $h=0$, we find
for any $1\leq q\leq \infty:$ 
\begin{eqnarray}
 \|p\|_{L^{\infty}(0,T;L^q(\Omega \times \mathbb{R}^N))} \!\leq\! 
 e^{[N\beta/q'+ \|Êa^- \|_\infty] T} & \hskip -4mm \Big[&
 \hskip -4mm \|p_0\|_{L^q(\Omega \times \mathbb{R}^N)} \!+\!
 \|g\|_{L^q_k(\Sigma^-_T)}  \Big], \hskip 4mm  \label{lp1} \\
\|{\rm Tr}\, p\|_{L^q_k(\Sigma^+_T)}  \!\leq\! 
  e^{[N\beta/q'+ \|Êa^- \|_\infty] T}& \hskip -4mm \Big[& 
  \hskip -4mm \|p_0\|_{L^q(\Omega \times \mathbb{R}^N)} \!+\!
 \|g\|_{L^q_k(\Sigma^-_T)}  \Big].  \hskip 4mm \label{lp2}
\end{eqnarray} 
% f^(p-1) h -> (f^p)1/p' (h^p)1/p
% 1/p + p-1/p = 1/p + 1/(p/(p-1))
% sale para q finito 
\end{itemize}
\end{theo}

The positivity result stated in Theorem \ref{positivity} implies a maximum principle.

\begin{theo}{\bf (Maximum principle).}
\label{comparison}
Under the hypotheses of Theorems \ref{positivity} and \ref{conservation}, the following two comparison principles hold:
\begin{itemize}
\item[(i)] if $p_1$ and $p_2$ are two solutions with data $h_1\leq h_2$, 
$g_1\leq g_2$, and $p_{1,0} \leq p_{2,0}$, then $p_1 \leq p_2$.
\item[(ii)] if $p_1$ and $p_2$ are two nonnegative solutions with the same data $h$, 
$g$, $p_0$, and coefficients $a_1=a_1^+ -a_1^-$, $a_2=a_1^+ -\|a_1^-\|_{\infty}$, 
so that $a_1^- \leq a_2^-$, then $p_1 \leq p_2$.
\end{itemize}
The results still hold true if ${\rm div}_{\mathbf v}( {\mathbf F} p)$ is replaced by ${\mathbf F} \cdot \nabla_{\mathbf v} p$, where ${\mathbf F}$ is a bounded field depending also on 
${\mathbf v},$ in such a way that ${\rm div}_{\mathbf v} \mathbf F$ is bounded. 
Moreover, if $g \in L^{\infty}(\Sigma^+_T)$ the solution $p$ satisfies:
\begin{eqnarray} 
\|p\|_{L^{\infty}(Q_T)} \leq e^{[N \beta + \|Êa^-\|_{\infty}]T}
\left[  \|p_0\|_{\infty} +  \|g\|_{\infty} + \int_0^t  \|h(s)\|_{\infty} ds\right]. \label{linf}
\end{eqnarray}
\end{theo}

{\bf Proof.} Let us first extend the positivity result in Theorem \ref{positivity}
to fields ${\mathbf F}$ depending on  ${\mathbf v}$.
We set $\overline{p}= e^{-(\lambda+N\beta) t} p({\mathbf x}, e^{-\beta t}{\mathbf v}, t)$ and $\overline{h}= e^{-(\lambda+N\beta) t} h({\mathbf x}, e^{-\beta t}{\mathbf v}, t)$. 
Then, $\overline{p}$ satisfies the equation:
\begin{eqnarray*}
\frac{\partial \overline{p}}{\partial t} \!+\! e^{-\beta t} \mathbf v\!\cdot\!\nabla_{\mathbf x}  \overline{p} \!+\! e^{\beta t}{\bf F}(\mathbf x, e^{-\beta t} \mathbf v, t) \!\cdot\! \nabla_{\mathbf v} \overline{p}
\!-\!\sigma e^{2 \beta  t}\Delta_{\mathbf v} \overline{p} \!+\! (a(\mathbf x, e^{-\beta t} \mathbf v, t) +\lambda) \overline{p} \!=\! \overline{h}.
\end{eqnarray*}
We multiply by $\overline{p}^-$ and integrate to get:
\begin{eqnarray*}
- \int_{\Omega \times {\mathbb R}^N} \hskip -5mm 
{|\overline{p}^-(T)|^2 \over 2} d \mathbf x d \mathbf v
- \int_{\partial \Omega \times {\mathbb R}^N \times [0,T]} \hskip -1.5cm
e^{-\beta t} {\mathbf v} \cdot {\mathbf n} {|\overline{p}^-|^2 \over 2}
d S d \mathbf v dt
+ \int_{\Omega \times {\mathbb R}^N \times [0,T]} \hskip -1.5cm
{\rm div}_{\mathbf v} {\bf F}(\mathbf x, e^{-\beta t} \mathbf v, t) 
{|\overline{p}^-|^2 \over 2} \nonumber \\
-\sigma  \int_{\Omega \times {\mathbb R}^N \times [0,T]} \hskip -1.2cm
e^{2 \beta  t} |\nabla_{\mathbf v} \overline{p}^-|^2 
-  \int_{\Omega \times {\mathbb R}^N \times [0,T]} \hskip -1.2cm
(a(\mathbf x, e^{-\beta t} \mathbf v, t)+\lambda) |\overline{p}^-|^2 =  
\int_{\Omega \times {\mathbb R}^N \times [0,T]} \hskip -1.2cm
\overline{h} \overline{p}^- \geq 0.
\end{eqnarray*}
Notice that $\overline{p}^-(0)=0$ for $p(0) \geq 0$ and $\overline{p}^-=0$ on 
$\Sigma^-_T$. The only contribution to the integral over $\partial \Omega$ comes
from the region where ${\mathbf v} \cdot {\mathbf n} >0$. Choosing 
$\lambda \geq  \|a^-\|_{\infty}+ \|{\rm div}_{\mathbf v} {\bf F} \|_{\infty}$, we conclude 
that $|\overline{p}^-|=0$. Therefore, $p \geq 0$ if $p_0 \geq 0$, $h \geq 0$ and
$p\big|_{\Sigma^-_T} \geq 0$.

Assertion (i) is a consequence of the positivity result. Indeed, setting $\overline{p}=p_2-p_1$, linearity plus the positivity result imply that $p_2 - p_1\geq 0.$ 

To prove statement (ii),  we set $\hat{p_1}=e^{-\|a_1^-\|_{\infty}t} p_1$ and $\hat{p_2}=e^{-\|a_1^-\|_{\infty}t} p_2$. These functions are solutions of similar problems, with source $\hat{h}=e^{-\|a_1^-\|_{\infty}t} h$, boundary datum $\hat{g}=e^{-\|a_1^-\|_{\infty}t} g$ and initial datum $p_0$:
\begin{eqnarray}
 \frac{\partial \hat{p}_1}{\partial t} \!+\!  {\mathbf v} \!\cdot\! \nabla_{\mathbf x} \hat{p}_1
 \!+\! {\bf F} \!\cdot\! \nabla_{\mathbf v} \hat{p}_1 \!-\! \beta {\rm div}_{\mathbf v}
 ( {\mathbf v}\hat{p}_1)  \!-\! \sigma\Delta_{\mathbf v} \hat{p}_1
 \!+\! a_1^+\hat{p}_1& \hskip -3mm = \hskip -3mm & 
 (a_1^- \!-\! \|a_1^-\|_{\infty}) \hat{p_1} \!+\! \hat{h}  \nonumber \\
 \frac{\partial \hat{p}_2}{\partial t} \!+\!  {\mathbf v} \!\cdot\! \nabla_{\mathbf x} 
 \hat{p}_2  
 \!+\! {\bf F} \!\cdot\! \nabla_{\mathbf v} \hat{p}_2 \!-\! \beta {\rm div}_{\mathbf v}
 ( {\mathbf v}\hat{p}_2) \!-\! \sigma\Delta_{\mathbf v} \hat{p}_2
 \!+\! a_1^+\hat{p}_2 &\ \hskip -3mm=\hskip -3mm & \hat{h}. \nonumber
\end{eqnarray}
Since $ (a_1^- -\|a_1^-\|_{\infty}) \hat{p}_1\leq 0$, assertion (i) implies that 
$p_1\leq p_2$.

For the $L^{\infty}$ estimate, let us first notice that  if $p$ is a solution with data 
$h, g, p_0 \leq 0$ then $-p$ is a solution with data $-h, -g, -p_0 \geq 0$ 
by linearity. Therefore, $-p \geq 0$ and $p \leq 0$.  The reverse inequality 
holds too. Now, let us set $p = e^{\lambda t} \hat{p}$ with 
$\lambda= N\beta + \|a^-Ê\|_{\infty}$. 
The function $\hat{p}$ is a solution 
of equations (\ref{lin:ec})-(\ref{lin:in}), (\ref{bc:plin})  with an additional source 
term $-\lambda e^{-\lambda t}  p= -\lambda  \hat{p}$. 
Set  $M(t) =  \int_0^t  e^{-\lambda s}\|h(s)\|_{\infty} ds
+  \|g\|_{\infty} + \|p_0\|_{\infty}$
and $\overline{p} = \hat{p}- M$. Then, $\overline{p}$ satisfies: 
\begin{eqnarray*}
 \frac{\partial \overline{p}}{\partial t} 
 + {\mathbf v} \!\cdot\! \nabla_{\mathbf x} \overline{p}+
 {\bf F} \!\cdot\! \nabla_{\mathbf v} \overline{p} 
 - \beta {\rm div}_{\mathbf v} ( {\mathbf v}\overline{p})  
 -\sigma\Delta_{\mathbf v} \overline{p}
 + (a + \lambda) \overline{p}  \nonumber \\
 %e^{-\lambda t}  h - e^{-\lambda t}  h  - (-N\beta + a^+ - a^-  + \|a^-\| + N \beta ) M 
 =  - a^+ M - (\|a^-\|_{\infty}- a^- ) M \leq 0
\end{eqnarray*}
with initial and boundary conditions $p_0-M \leq 0$ and $e^{-\lambda t}g-M \leq 0$. Notice that $e^{-\lambda t}<1$ because $-\lambda <0$.
Therefore, $\overline{p} \leq 0$, $\hat{p} \leq M$ and $p \leq e^{\lambda t} M$. The reverse inequality follows in a similar way by linearity.

\section{Estimates on velocity integrals}
\label{sec:decayv}

The nonlinear problem includes the velocity integrals $\rho(p)$ and $j(p)$ of the density $p$ as coefficients. In this section we discuss strategies to estimate
velocity integrals in terms of density norms. Let us start with the variable $j$.

\begin{lemma} \label{cotaj}
For any  $p \geq 0$, the norms 
$\| j \|_{L^q_{\mathbf x}}$, $1 \leq q \leq \infty,$ 
of the flux $j$ defined in equation (\ref{eq:intpintvp}) can be bounded
in terms of $\|p \|_{L^\infty_{\mathbf x \mathbf v}}$.
\end{lemma}

{\bf Proof.} 
Let us set $|\mathbf v| w(\mathbf v) = |\mathbf v| 
[1 + e^{|\mathbf{v}- \mathbf{v}_0 \chi |^2 / \sigma_v^2} ]^{-1}$. 
This function is bounded and integrable. 
Then,
\begin{eqnarray}
\| j \|_{L^\infty_{\mathbf x}} \leq \||\mathbf v|w\|_{L^1_{\mathbf v}}
\|p \|_{L^\infty_{\mathbf x \mathbf v}}, \label{jinfty} \\
\|j \|_{L^q_{\mathbf x}} \leq {meas}(\Omega)^{1/q} 
\|j\|_{L^\infty_{\mathbf x \mathbf v}}, \quad 1 \leq  q < \infty. \label{jq}  \\
\end{eqnarray}

The anastomosis term may be controlled using the kinetic equation, 
as we show below.

\begin{lemma} Under the hypotheses of Theorems \ref{positivity} and \ref{conservation}, let $p$ be a nonnegative solution of  
problem (\ref{eq:pnu}),(\ref{eq:pc0}),(\ref{eq:alphaF}) with boundary 
condition (\ref{bc:plin}) and nonnegative data. Assume that
$c \geq 0$. Then,  $\| \int_0^T \int p d \mathbf v ds \|_{L^2_{\mathbf x}}$ is bounded by the parameters of the problem and $\|p\|_{L^{\infty}_t L^\infty_{\mathbf x \mathbf v} }.$
\end{lemma}

{\bf Proof. }
Let us recall the equation of mass conservation from Theorem \ref{positivity}:
\begin{eqnarray*}
{\partial  \over \partial t}  \int \int p d {\mathbf v} d {\mathbf x}
+ \gamma \int \left[  \int_0^{t} \int p d {\mathbf v}'  dt'  \right] \left[\int
p  d {\mathbf v}\right] d {\mathbf x} = \int\int \alpha(c) \nu p d {\mathbf v} d {\mathbf x} 
\nonumber \\
+ \int_{\Gamma^-} \!\!\!\! |{\mathbf v} \cdot {\mathbf n}({\mathbf x}) |
g \,dS\, d{\mathbf v} 
- \!\! \int_{\Gamma^+} \!\!\!\! |{\mathbf v} \cdot {\mathbf n}({\mathbf x}) |
{\rm Tr} \,p \,dS\, d{\mathbf v}.
\label{massconservation2}
\end{eqnarray*}
Setting $a({\mathbf x},t)=\int_0^{t} \int p({\mathbf x},{\mathbf v}',t',) d {\mathbf v}'  dt' $, we notice that $ {da \over dt}({\mathbf x},t)= 
\int p({\mathbf x},{\mathbf v}',t)  d {\mathbf v}' $. Therefore:
\begin{eqnarray*}
\left[  \int_0^{t} \int p d {\mathbf v}'  dt'  \right] 
\int p({\mathbf x},{\mathbf v},t) d {\mathbf v}
= a({\mathbf x},t){da \over dt}({\mathbf x},t) = {1\over 2} 
{da^2 \over dt}({\mathbf x},t). \label{da}
\end{eqnarray*}
Integrating (\ref{massconservation2}) in time and inserting (\ref{da}), we find:
\begin{eqnarray*}
\int \int p (t) d {\mathbf v} d {\mathbf x} - \int \int p (0) d {\mathbf v} d {\mathbf x}
+ {\gamma \over 2} \int a({\mathbf x},t)^2 d {\mathbf x}-
 {\gamma \over 2} \int a({\mathbf x},0)^2 d {\mathbf x} =
\nonumber \\
\int_0^t \int \int  \alpha(c) \nu p \, ds d {\mathbf v} d {\mathbf x}
+ \int_{\Sigma^-} \!\!\!\! |{\mathbf v} \cdot {\mathbf n}({\mathbf x}) |
g \,dS\, d{\mathbf v} 
- \!\! \int_{\Sigma^+} \!\!\!\! |{\mathbf v} \cdot {\mathbf n}({\mathbf x}) |
{\rm Tr} \,p \,dS\, d{\mathbf v}.
\label{massconservation3}
\end{eqnarray*}
Notice that $a(0,{\mathbf x})^2=0$. Therefore:
\begin{eqnarray*}
\int \!\! d {\mathbf x} \left[ \int_0^t \!\! \int p d {\mathbf v}  ds \right]^2 \leq 
C(\gamma,\alpha_1,{\rm meas}(\Omega),\|p\|_{\infty}, \|\nu\|_{L^\infty(0,T,L^1_{\mathbf x \mathbf v})},
\|g\|_{L^1(\Sigma_T^-)}).
\label{cotarho2}
\end{eqnarray*}

To ensure the positivity of the solutions of linearized versions of equation 
(\ref{eq:pnu}), $b(p)=\int_0^t \rho(p) ds$ is taken to be a known coefficient.
To apply Theorem \ref{positivity}, it should be a bounded function. 
$L^q_{\mathbf x}$ estimates of $\rho(p)$ are obtained
controlling the moments.
%However, the proof of
%Theorem \ref{comparison} might yield positivity under the condition that
%$\int_{\Omega \times \mathbb R^N \times (0,T)} b(p_n) p^2 d \mathbf x
%d \mathbf v dt$ is finite but you must know that b_n is positive.
% and for existence you need boundedness anyway

\begin{lemma}
\label{moments}
Under the hypotheses of Theorems \ref{positivity} and \ref{conservation}, 
let $p$ 
a nonnegative solution of the linear equations (\ref{lin:ec})-(\ref{lin:in}) with boundary condition (\ref{bc:plin}) and nonnegative data.
If $(1+|\mathbf v|^2)^{\mu/2} p_0 \in L^1(\Omega \times \mathbb R^N)$,
$(1+|\mathbf v|^2)^{\mu/2} h \in L^\infty(0,T;L^1(\Omega \times \mathbb R^N))$
and $(1+|\mathbf v|^2)^{\mu/2} g \in L^1_k(\Sigma_T^-)$ for a positive integer
$\mu$, then, for $\ell=0,1,\ldots,\mu$ and $t \in [0,T]$, all the moments 
\begin{eqnarray*}
m^\ell(p(t))=  \int_ {\Omega \times \mathbb{R}^N}  \hskip -4mm
 |{\mathbf v}|^\ell p \,d{\mathbf x}\, d{\mathbf v}, \quad
m^\ell_k({\rm Tr} \, p^+)=  \int_ {\Sigma_T^+} 
|\mathbf v \cdot \mathbf n| |{\mathbf v}|^\ell {\rm Tr} 
\, p^+ \,dS\, d{\mathbf v} \, dt,
\end{eqnarray*} are bounded
in terms of the parameters $\beta$, $\sigma$, $N$, $T$, $\mu$, the norms of the 
data
$\|(1+|\mathbf v|^2)^{\mu/2}p_0\|_{L^1_{\mathbf x \mathbf v}}$,
$\|(1+|\mathbf v|^2)^{\mu/2}h\|_{L^\infty(0,T; L^1_{\mathbf x \mathbf v})}$,
$\|(1+|\mathbf v|^2)^{\mu/2} g\|_ {L^1_k(\Sigma_T^-)}$,
$\|a^-\|_{\infty}$, $\|Ê\mathbf F\|_{\infty},$
and  $\|p\|_{L^\infty(0,T;L^\infty_{\mathbf x \mathbf v})}$,
$\|p\|_{L^\infty(0,T;L^1_{\mathbf x \mathbf v})}.$
\end{lemma}

{\bf Proof.} Notice that the integral
\[
\int_ {\mathbb{R}^N}  { p \over |{\mathbf v}| } 
d{\mathbf v}  \leq  \|p\|_{L^\infty_{\mathbf x \mathbf v}}
\int_ {|{\mathbf v}| < R}  { d \mathbf v \over |{\mathbf v}| } 
+ {1 \over R}
\int_ {\mathbb R^N} \hskip -4mm  p d{\mathbf v}
\leq 
\|p\|_{L^\infty_{\mathbf x \mathbf v}}
{R^{N-1} \over N-1}
+ {1 \over R}
\int_ {\mathbb R^N} \hskip -4mm  p d{\mathbf v}.
\]
Since $\Omega$ is a bounded set, 
$\int_ {\Omega \times \mathbb{R}^N}  { p \over |{\mathbf v}| } 
d{\mathbf x} d{\mathbf v}$  is bounded in terms of 
$\|p\|_{L^\infty_{\mathbf x \mathbf v}}$ and
$\|p\|_{L^1_{\mathbf x \mathbf v}}$.
We first apply identity (\ref{momentum}) with $\mu=1$  to find:
\begin{eqnarray}
{d\over dt} \int_ {\Omega \times \mathbb{R}^N} \!\!\!\!\!\!\!\!\!\! 
|{\mathbf v}| p \,d{\mathbf x}\, d{\mathbf v}  = 
 \!\int_{\Gamma^-} \!\!\! |{\mathbf v} \cdot {\mathbf n}({\mathbf x}) | 
 |{\mathbf v}| g \,dS\, d{\mathbf v}
+\int_ {\Omega \times \mathbb{R}^N}   \!\!\!\!\!\!\!\! |{\mathbf v}| h \,d{\mathbf x}\, d{\mathbf v} \nonumber \\
 - \!\int_{\Gamma^+}\!\!\!\!\! |{\mathbf v} \!\cdot\!  {\mathbf n}({\mathbf x}) | 
 |{\mathbf v}| {\rm Tr} \,p \,dS\, d{\mathbf v} 
 -\beta \!\int_ {\Omega \times \mathbb{R}^N}\!\!\!\!\!\!\!\!\!\!  
 |{\mathbf v}|  p \,d{\mathbf x}\, d{\mathbf v}  
 -  \!\int_ {\Omega \times \mathbb{R}^N}\!\!\!\!\!\!\!\!\!\! 
a  |{\mathbf v}| p \,d{\mathbf x}\, d{\mathbf v} \nonumber \\
+ (N-1) \sigma  \!\!\int_ {\Omega \times \mathbb{R}^N}\!\!\!\!\!\!\!\!\!\! 
|{\mathbf v}|^{-1}  p \,d{\mathbf x}\, d{\mathbf v}
 +  \!\int_ {\Omega \times \mathbb{R}^N}\!\!\!\!\!\!\!\!\!\!  
{\mathbf F} \cdot {\mathbf v} |{\mathbf v}|^{-1}  p \,d{\mathbf x}\, d{\mathbf v}.
\nonumber
\end{eqnarray}
Integrating in time,  we find:
\[
\int_ {\Omega \times \mathbb{R}^N} \!\!\!\!\!\!\!\!\!\! 
|{\mathbf v}| p \,d{\mathbf x} d{\mathbf v} ds \leq
C(p_0,g,h,\mathbf F,p) + \|a^-\|
\int_0^t \int_ {\Omega \times \mathbb{R}^N} \!\!\!\!\!\!\!\!\!\! 
|{\mathbf v}| p \,d{\mathbf x} d{\mathbf v} ds, 
\]
where $C(p_0,g,h,\mathbf F,p)$ depends on the norms and parameters
mentioned in the statement. Gronwall's Lemma provides the required bound 
on $\int_ {\Omega \times \mathbb{R}^N}  \!\!
|{\mathbf v}| p \,d{\mathbf x} d{\mathbf v} ds$. Once the moment of $p$ is 
bounded,
the estimate on the moment of its trace follows inserting this information in the
differential equation.

We reason by induction. Assuming that the moments $m_\ell(p)$
are bounded in terms of the required norms for $\ell \leq M-1$,
let us see that the same holds true of $m_M(p)$.
Integrating in time (\ref{momentum}), using the bounds
on $m_{M-1}(p)$ and $m_{M-2}(p)$, together with Growall's
lemma, we find the desired estimate. By induction, it holds
up to $M=\mu$. Once the moments of $p$ are bounded,
the estimate on the moment of its trace follows inserting this 
information in the differential equations.
\\

The relation between velocity moments and norms of the
marginal density $\rho(p)= \int_{\mathbb R^N} p d \mathbf v$ is 
established in the following lemma
 \cite{unbounded,degond}.
All $L^q_{\mathbf x}$ norms for finite $q$ can be controlled
in that way. To obtain $L^\infty_{\mathbf x}$ estimates of the marginal 
density $\rho(p)$ we resort to a strategy 
involving velocity weights introduced in reference \cite{degond}.

\begin{lemma} 
\label{interp} Let $\Omega \subset \mathbb R^N$ be bounded.
For any nonnegative $p$ the following inequalities hold:
\begin{eqnarray}
 \| |\mathbf v|^{\ell} p \|_{L^1(\Omega \times \mathbb R^N)}  \leq 
\| p \|_{L^1(\Omega \times \mathbb R^N)}^{1- {\ell \over \mu}} \;
\| |\mathbf v|^\mu p \|_{L^1(\Omega \times \mathbb R^N)}^{{\ell \over \mu}}, 
\quad \mu > \ell >0, \label{boundmlp1} \\
\| \int_{I \!\! R^N} \hskip -3mm |\mathbf v|^\ell p \, d {\mathbf v} 
\|_{L^{N+\mu \over N+\ell}(\Omega)}  \leq C_{N,\mu,\ell} \,
\| p \|_{L^\infty(\Omega \times \mathbb R^N)}^{\mu-\ell \over N+\mu} \;
\| |\mathbf v|^{\mu} p \|_{L^1(\Omega \times \mathbb R^N)}^{N+\ell \over N+\mu}, 
\quad \mu > \ell >0,  \label{boundmlpinf} \\
\|\int_{{\mathbb R}^N} \hskip -3mm 
|{\mathbf v}| p d{\mathbf v} \|_{L^\infty(\Omega)} \leq  
C_\mu \|p\|_{L^\infty(\Omega \times \mathbb R^N)}^{1-(N+1)/\mu}
\| (1\!+\!|{\mathbf v}|^2)^{\mu\over 2} p\|_{L^\infty(\Omega \times 
\mathbb R^N)}^{(N+1)/\mu}, 
\quad \mu>N+1,
\label{boundvinf} \\
\| \int_{{\mathbb R}^N} \hskip -3mm p d {\mathbf v} \|_{L^\infty(\Omega)}
\leq  C_\mu  \| p \|_{L^\infty(\Omega \times \mathbb R^N)}^{1-N/\mu}
\|(1\!+\!|{\mathbf  v}|^2)^{\mu\over 2} p \|_{L^\infty(\Omega \times 
\mathbb R^N)}^{N/\mu}, 
\quad \mu>N,
\label{boundinf} \\
\|(1\!+\!|{\mathbf  v}|^2)^{\mu-1 \over2} \! p \|_{L^\infty(\Omega \times \mathbb R^N\!)} 
\leq C_\mu \|p\|_{L^\infty(\Omega \times \mathbb R^N\!)}^{1/\mu} 
\|(1\!+\!|{\mathbf  v}|^2)^{\mu\over 2} p \|_{L^\infty(\Omega \times \mathbb R^N\!)}^{1-1/\mu}, 
\; \mu>1, 
\label{boundinterp}
\end{eqnarray}
provided the involved integrals and norms are finite.
\end{lemma}

Revising the proof of this lemma in reference \cite{unbounded}, we see
that it extends to the traces on the boundary, with respect to either the
Lebesgue or the kinetic measure.

\begin{cor}
\label{interptrace}
The inequalities in Lemma \ref{interp} hold for ${\rm Tr} \,p^+$
replacing the spaces
$L^q_{\mathbf x}L^1_{\mathbf v}(\Omega \times \mathbb R^N)$ 
by  $L^q_{\mathbf x}L^1_{\mathbf v}(\Gamma^+)$
and $L^q(\Omega \times \mathbb R^N)$ by $L^q(\Gamma^+)$
provided the involved integrals and norms are finite.
\end{cor}

\begin{cor}
\label{interptracekin}
The inequalities in Lemma \ref{interp} hold for $|\mathbf v \cdot
\mathbf n|{\rm Tr} \,p^+$ replacing the spaces
$L^q_{\mathbf x}L^1_{\mathbf v}(\Omega \times \mathbb R^N)$ 
by  $L^q_{\mathbf x}L^1_{\mathbf v}(\Gamma^+)$
and $L^q(\Omega \times \mathbb R^N)$ by $L^q(\Gamma^+)$
provided the involved integrals and norms are finite.
\end{cor}

%
% We use  the L_v estimate for p \chi_{\Gamma^+}
% then integrate in \partial \Omega
% 
% We apply it to |v n| Tr p^+
%

Let us now estimate the velocity decay of $p$,
and the $L^\infty_{\mathbf x}$ norms of velocity integrals, which
extend to traces on the boundary.

\begin{prop} 
\label{velocityLinf}
Let $p \geq 0$ be a solution of the initial value problem (\ref{lin:ec})-(\ref{lin:in})
with  boundary conditions given by (\ref{bc:plin}).
Under the hypotheses:
\begin{itemize}
\item [(i)] $a \in L^{\infty}(\Omega \times I\!\!R^N \times (0,T))$, 
\item[(ii)] $(1+|{\mathbf v}|^2)^{\mu/2} p_0({\mathbf x},{\mathbf v}) \in
L^1\cap L^{\infty}(\Omega \times {\mathbb R}^N),  \; \mu >N, 
\quad p_0 \geq 0$, 
\item[(iii)] $(1+|{\mathbf v}|^2)^{\mu/2} 
g({\mathbf x},{\mathbf v},t) \in L^1 \cap L^{\infty}(\Sigma^-_T),$ 
%\item[(iv)]
%$(1+|{\mathbf v}|^2)^{\mu/2} h({\mathbf x},{\mathbf v},t) \in
%L^1 \cap L^{\infty}(Q_T),$
\item[(iv)] $\mathbf F \in L^{\infty}(\Omega \times (0,T))$,
\end{itemize}
the norms
$\| (1+|{\mathbf v}|^{2})^{\mu/2} p \|_{L^\infty(0,T;L^{\infty}_{\mathbf x \mathbf v})
}$,  
and $\| p\|_{L^\infty(0,T;L^\infty_{\mathbf x}L^1_{\mathbf v})}$ 
are bounded by constants depending on  
$T$, $\sigma$, $\beta$, $\mu$, as well as 
$\| (1+|{\mathbf v}|^{2})^{\mu/2} p_0 \|_{L^{\infty}_{\mathbf x \mathbf v }}$, 
%$\| (1+|{\mathbf v}|^{2})^{\mu/2} h \|_{L^{\infty}(Q_T)}$,
$\| (1+|{\mathbf v}|^{2})^{\mu/2}  g \|_{L^\infty(\Sigma^-_T)}$, 
$\|a^-\|_{L^\infty(Q_T)}$, $\|{\mathbf F}\|_{L^{\infty}(\Omega \times (0,T))}$
and $\| p\|_{L^\infty(0,T;L^\infty_{\mathbf x \mathbf v})}.$
%If $\mu > N+1$, the norm $\||{\mathbf v}| p\|_{L^\infty(0,T;L^\infty_{\mathbf x}L^1_{\mathbf v})}$ is similarly bounded. 
Moreover, if
\begin{itemize}
\item[(v)] $(1+|{\mathbf v}|^2)^{\mu/2} 
g({\mathbf x},{\mathbf v},t) \in  L^1_k \cap L^{\infty}_k(\Sigma^-_T),$ 
\end{itemize}
then, for any $1\leq q\leq \infty:$ 
\begin{eqnarray}
 \|(1+|{\mathbf v}|^2)^{\mu\over 2} p\|_{L^{\infty}(0,T;L^q_{\mathbf x \mathbf v})} 
 \!\leq\! e^{ \|ÊD \|_\infty T}  \Big[
  \|p_0\|_{L^q_{\mathbf x \mathbf v}} \!+\!
 \|(1+|{\mathbf v}|^2)^{\mu\over 2} g\|_{L^q_k(\Sigma^-_T)}  \Big],  
  \label{lp1weight} \\
\|(1+|{\mathbf v}|^2)^{\mu\over 2}  {\rm Tr}\, p\|_{L^q_k(\Sigma^+_T)}  \!\leq\! 
  e^{ \|ÊD \|_\infty T}  \Big[ 
  \|p_0\|_{L^q(\Omega \times \mathbb{R}^N)} \!+\!
 \|(1+|{\mathbf v}|^2)^{\mu\over 2}  g\|_{L^q_k(\Sigma^-_T)}  \Big],   
 \label{lp2weight}
\end{eqnarray} 
where $\|ÊD \|_\infty $ depends on $\sigma$, $\beta$, $\mu$,  $N$,
$\|a^-\|_{L^\infty(Q_T)}$ and $\|{\mathbf F}\|_{L^{\infty}(\Omega \times (0,T))}$.

\end{prop}

% for $t \in [0,T]$:
%\begin{eqnarray}
%\| Y(t) \|_{L^{\infty}_{\mathbf x \mathbf v}} \leq 
%M(T,\sigma,\gamma,k,\eta,d_1,p_0,c_0,f,\|a^-\|_{L^\infty_{\mathbf x \mathbf v}},
%\|p\|_{L^\infty_{\mathbf x \mathbf v}}, \|p\|_{L^1_{\mathbf x \mathbf v}},
%\| |{\mathbf v}|^2 p\|_{L^1_{\mathbf x \mathbf v}}). \nonumber
%\end{eqnarray}}

{\bf Proof.} We set $Y({\mathbf x},{\mathbf v},t)=(1+|{\mathbf v}|^2)^{\mu/2} 
p({\mathbf x},{\mathbf v},t)$.
Multiplying  equation (\ref{lin:ec}) by $(1+|{\mathbf v}|^2)^{\mu/2}$, 
$\mu>0$, we get:
\begin{eqnarray}
{\partial \over \partial t} Y \!+\! {\mathbf v} \nabla_{\mathbf x} Y \!+\! \Big( {\mathbf F} 
+ 2 \sigma \mu {{\mathbf v} \over 1 + |{\mathbf v}|^2} \!-\! \beta {\mathbf v} \Big) 
\nabla_{\mathbf v} Y \!-\! \Delta_{\mathbf v} Y \!=\! (N \beta \!-\!a) Y 
\!+\! R \label{eq:weight}
\end{eqnarray}
where $R=R_1+R_2+R_3$, with
\begin{eqnarray}
R_1 = %(1+|{\mathbf  v}|^2)^{\mu/2} h + 
{\mu} (1+|{\mathbf v}|^2)^{\mu/2 -1}  {\mathbf F} \cdot {\mathbf v} p,  \quad
R_2 =  -\beta \mu {|{\mathbf  v}|^2 \over (1+|{\mathbf v}|^2)} Y,  \nonumber \\
R_3= \sigma \mu (\mu+2)   {|{\mathbf v}|^2 \over (1+|{\mathbf v}|^2)^2} Y - 
N \sigma \mu{1 \over 1+|{\mathbf v}|^2 } Y.
\nonumber
\end{eqnarray}
Thanks to Theorem \ref{comparison}:
\begin{eqnarray}
\| Y(t) \|_{L^\infty_{\mathbf x \mathbf v}} 
\!\!\leq C(p_0,g)
\!\!+\!\!\! \int_0^t \hskip -2mm\Big[ \! 
[N \beta \!+\!\|a^-\|_{\infty}  ] \| Y \|_{L^\infty_{\mathbf x \mathbf v}} 
\!\!+\!\! \|R_1 \|_{L^\infty_{\mathbf x \mathbf v}} 
\!\!+\!\! \|R_2 \|_{L^\infty_{\mathbf x \mathbf v}} 
\!\!+\!\! \|R_3 \|_{L^\infty_{\mathbf x \mathbf v}} \!\! \Big] ds, \nonumber
\end{eqnarray}
% The proof in Degond requires data in Linf but also L2. Putting L1
% takes care of it.
where $C(p_0,g)$ is a constant depending on $\| (1+|{\mathbf  v}|^2)^{\mu/2} p_0\|_{\infty}$ and $\| (1+|{\mathbf  v}|^2)^{\mu/2} g\|_{\infty}$.
The factors ${|{\mathbf v}|^\varepsilon \over 1+|{\mathbf v}|^2} \leq 1$, for
$0 \leq \varepsilon \leq 2$. Therefore,
\begin{eqnarray}
 \|R_1\|_{L^\infty_{\mathbf x \mathbf v}} \leq  
 %\|(1+|{\mathbf  v}|^2)^{\beta/2} h \|_{L^\infty_{\mathbf x \mathbf v}}  +
 \mu \|(1+|{\mathbf v}|^2)^{\mu/2 -1}  {\mathbf F} \cdot {\mathbf v} p 
 \|_{L^\infty_{\mathbf x \mathbf v}}, \nonumber \\
 \|R_2\|_{L^\infty_{\mathbf x \mathbf v}} \leq   
  \beta \mu \|Y \|_{L^\infty_{\mathbf x \mathbf v}}, \quad
 \|R_3\|_{L^\infty_{\mathbf x \mathbf v}} \leq   \sigma \mu (\mu + 2 + N) 
 \|Y \|_{L^\infty_{\mathbf x \mathbf v}}. \nonumber
\end{eqnarray}
To bound
$\|(1+|{\mathbf v}|^2)^{\mu/2 -1}  {\mathbf F} \cdot {\mathbf v} p 
\|_{L^\infty_{\mathbf x \mathbf v}}$, we set:
\begin{eqnarray*}
\|(1+|{\mathbf v}|^2)^{\mu/2 -1}  
{\mathbf F} \cdot {\mathbf v} p \|_{L^\infty_{\mathbf x \mathbf v}}
\leq  { N |{\mathbf v}| \over 1+|{\mathbf v}|^2}
\|{\bf F}\|_{\infty} \|Y\|_{L^\infty_{\mathbf x \mathbf v}}
\leq N \|{\bf F}\|_{\infty} \|Y\|_{L^\infty_{\mathbf x \mathbf v}}.
\label{boundF0}
\end{eqnarray*}
Taking 
$A= (N \|{\bf F}\|_{\infty}+\beta) \mu 
+ \sigma \mu (\mu + 2 + N) + N \beta + \|a^-\|_{\infty} $
%+ \|(1+|{\mathbf  v}|^2)^{\beta/2} h \|_{L^\infty_{\mathbf x \mathbf v}} $ 
and  $B = C(p_0,g)$, 
%+ \|(1+|{\mathbf  v}|^2)^{\beta/2} f \|_{L^{1}(0,T; 
%L^{\infty}_{\mathbf x \mathbf v})}$
Gronwall's inequality implies
\begin{eqnarray}
\| Y(t) \|_{L^\infty_{\mathbf x \mathbf v}} \leq B e^{A t}, \quad t \in [0,T].
\nonumber 
\end{eqnarray}

Once the velocity decay has been established, the $L^{\infty}$ bounds on 
$ \int_{\mathbb R^N} p d {\mathbf v}$ follow from inequality
(\ref{boundinf}) in Lemma \ref{interp}. 
%When $\beta > N+1$, the estimate
%on $ \int_{\mathbb R^N} |\mathbf v| p d {\mathbf v}$ follows from
%inequality (\ref{boundvinf}). 

Writing down the analogous of equation (\ref{lp}) for equation (\ref{eq:weight})
we find:
\begin{eqnarray*} \hskip -8mm
\frac{d}{d t}\|Y(t)\|^q_{L^q(\Omega\times \mathbb{R}^N)} \hskip -3mm &=&
\hskip -3mm
\int_{\Gamma_{-}}|{\mathbf v}\cdot {\mathbf n}({\mathbf x})|
(1+ |\mathbf v|^2)^{\mu/2} g^q dS \, d{\mathbf v}  \nonumber \\
&-&  \hskip -3mm \int_{\Gamma_{+}}|\mathbf v\cdot \mathbf n(\mathbf x)|
(1+ |\mathbf v|^2)^{\mu/2} {(\rm Tr} \,p)^q  dS \, d{\mathbf v} 
\nonumber \\
&-& \sigma q(q\!-\!1) \!\!\int_{\Omega\times \mathbb{R}^N} \hskip -6mm 
Y^{(q\!-\!2)}|\nabla_{\mathbf v} Y|^2  d{\mathbf x} \, d{\mathbf v}
-q \int_{\Omega \times \mathbb R^N} 
\hskip -4mm D \, Y^q d{\mathbf x} \, d{\mathbf v},
\end{eqnarray*}
where $D$ is a bounded coefficient. Integrating in time, we
recover estimates (\ref{lp1}) and (\ref{lp2}) for $Y$ updating the
data, and replacing the exponent of the exponential by 
$\|DÊ\|_{\infty}.$

\section{Coupling to the diffusion equation with Neumann boundary condition}
\label{sec:heat}

In this section, we consider diffusion problems of the form:
\begin{eqnarray}
\frac{\partial}{\partial t}c(\mathbf{x},t) = d  \Delta_{\mathbf x} c(\mathbf{x},t) - \eta c(\mathbf{x},t)j(\mathbf{x},t) + h(\mathbf{x},t) ,  & {\mathbf x} \in \Omega, t>0, 
\label{eq:ch}\\
{\partial c \over \partial r}(\mathbf x,t) = c_{r_0}(\mathbf x,t), \; \mathbf x \in S_{r_0},
\quad {\partial c \over \partial r}(\mathbf x,t) = 0, \; \mathbf x \in S_{r_1}, &   t>0,
\label{eq:pc0h} \\
c(\mathbf{x},0) =c_0(\mathbf{x}),  &    {\mathbf x} \in \Omega,
\label{bc:ch}
\end{eqnarray}
where $d , \eta >0$, $c_{r_0}<0$ and
$ j(\mathbf{x},t)={j}(p)= \int_{\mathbb R^2} {|\mathbf{v}| 
\over 1 + e^{|\mathbf{v}- \mathbf{v}_0 \chi |^2/\sigma_v^2}}
p(\mathbf{x},\mathbf{v},t)\, d \mathbf{v} $. 
The domain $\Omega = \{ \mathbf x \in  \mathbb{R}^N \, | \, 
r_0 < r=|{\mathbf x}| < r_1 \}$, with boundaries
$S_{r_0} = \{ \mathbf x \in  \mathbb{R}^N \, \big| \, 
|\mathbf x| =r_0 \}$ and $S_{r_1} = \{ \mathbf x \in  \mathbb{R}^N \, \big| \, 
|\mathbf x| =r_1 \}.$

When $j(\mathbf{x},t)$ is a bounded function, existence of a unique 
global solution for equations (\ref{eq:ch})-(\ref{bc:ch}) can be proved by classical galerkin or spectral methods \cite{galerkin}. Coercitivity of the associated bilinear form is not necessary. However, it holds whenever $j(t)$ is continuous and does not vanish identically for any $t\in [0,T]$.

Let us now establish comparison and maximum principles that will be essential in the sequel.

\begin{prop} {(\bf Comparison principle).}
\label{comparisonneuman}
Let $c\in C([0,T];L^2(\Omega))$ be a solution of problem (\ref{eq:ch})-(\ref{bc:ch}) with initial datum $c_0\in L^2(\Omega)$, boundary condition $c_{r_0} \in C([0,T];L^2(\partial \Omega))$ and nonnegative coefficient $j\in 
L^{\infty}(\Omega \times (0,T))$. If $c_0\geq 0$, $h \geq 0$ and $c_{r_0} \leq 0$, 
then $c\geq 0$. Moreover, the following comparison principle holds.
Given two solutions $c_1$ and $c_2$ with sources $h_1,h_2$,
initial data $c_{1,0}, c_{2,0}$ and normal derivatives at the boundary $g_1,g_2$,
if $g_1 \leq  g_2$, $c_{1,0} \leq c_{2,0}$,  $h_1 \leq h_2$, 
then $c_1 \leq c_2$.

\end{prop}

{\bf Proof.} 

Multiplying the equation 
\begin{eqnarray}
\frac{\partial}{\partial t} {c}(\mathbf{x},t) = d \Delta_{\mathbf x} 
 {c}(\mathbf{x},t) - \eta  {c}(\mathbf{x},t) j(\mathbf{x},t) + {h},
\nonumber
\end{eqnarray}
by $ {c}^-={\rm Max}(-c,0)$ and integrating, we get
\begin{eqnarray}
{1\over 2} \| {c}^-(t)\|_2^2 + \!\! \int_0^t \!\!\int_{\Omega} [ |\nabla {c}^-|^2  
 +  \eta j  | {c}^-|^2 ]  = \nonumber \\
{1\over 2} \| {c}^-(0)\|_2^2 
- \!\! \int_0^t \!\! \int_{\partial \Omega} {\partial  {c} \over \partial 
{\mathbf n}} 
c^{-} \!\! - \!\! \int_0^t \!\! \int_{\Omega} h c^- \leq 0,  \label{cneg}
\end{eqnarray}
% c_x n_x c = 1/2 (c^2)_x
since, in our case, 
\begin{eqnarray}
 - \int_{\partial \Omega} {\partial c \over \partial {\mathbf n}} c^{-}
 =   - \int_{r=r_1} {\partial c \over \partial r}(r_1) c^{-}
 + \int_{r=r_0}  {\partial c \over \partial r}(r_0) c^{-} =  
 \int_{r=r_0}  {\partial c \over \partial r}(r_0) c^{-} \leq 0.
\nonumber
\end{eqnarray}
This implies that $c^-=0$ and $c\geq 0$.

If $h \geq 0$, ${\partial c \over \partial {\mathbf n}} \geq 0$ and $c(0) \geq 0$,
inequality (\ref{cneg}) implies immediately $c \geq 0.$ Reproducing the
computations for $\overline{c}= {c_2}- {c_1}$, inequality (\ref{cneg}) 
implies $c_2 \geq c_1$ by linearity. 

\begin{cor}
If $c$ is a solution of equations (\ref{eq:c})-(\ref{eq:pc0}),(\ref{bc:c}) with nonnegative data $c_0$ and coefficient $j$, then $c\geq 0$ and $c \leq u$, 
$u$ being the solution of the heat equation with the same initial and 
boundary data, but zero source.
\end{cor}
{\bf Proof.}
Positivity is a straightforward consequence of the previous maximum principle.
Similarly, the comparison principle applied with $h_1=-c j$ and 
$h_2=0$ implies $c \leq u$. \\

To control the tumor angiogenic factor (TAF) induced force field
$\|{\mathbf F}(c)\|_{L^{\infty}_{\mathbf x t}}$, we will need $L^{r}-L^{q}$ estimates 
of $c$ analogous to the known estimates for solutions of heat equations in the 
whole space. Let us first consider the pure initial value problem:
\begin{eqnarray}
u_t(\mathbf x,t)= d \Delta_{\mathbf x} u(\mathbf x,t),  
&  {\mathbf x} \in \Omega, \, t>0,\label{ec:uin} \\
{\partial u \over \partial {\bf n}}(\mathbf x,t) =  0, &   \,
\mathbf x \in S_{r_0}\cup S_{r_1}, \,  t>0, \label{ec:bcuin} \\
u(\mathbf{x},0) =u_0(\mathbf{x}),  & \quad  {\mathbf x} \in \Omega.
\label{ec:inuin}
\end{eqnarray}
For any $u_0 \in L^2(\Omega)$, there is a unique global solution 
$u \in C([0,T],L^2(\Omega)) \cap L^2(0,T; H^1(\Omega)),$ see
reference \cite{galerkin}. $\Omega$ being bounded, this
remains true when $u_0 \in L^\infty(\Omega).$
We can construct the solution using eigenfunction expansions.
Let $\phi_n$, $n=1,2...$, be the orthonormalized eigenfunctions for 
the homogeneous Neumann problem:
\begin{eqnarray} 
- d \Delta \phi_n = \lambda_n \phi_n \quad \mbox {on} \; \Omega, \quad
\frac{\partial}{\partial {\bf n}} \phi_n=0 \quad \mbox{on} \; \partial \Omega.
\label{eigenvalue}
\end{eqnarray}
The smallest eigenvalue is $\lambda_1=0$ with constant eigenfunction.
Using separation of variables, $u$ takes the form:
\begin{eqnarray}
u(\mathbf{x},t) = \sum_{n\geq 1} u_{0,n} \phi_n(\mathbf{x}) e^{-\lambda_n t},
%+  \sum_{n\geq 1} \phi_n(\mathbf{x}) \int_0^t \tilde{z}_n(s) e^{-\lambda_n (t-s)} ds, \nonumber
\quad
u_{0,n}  = \int_{\Omega}   u_0(\mathbf{y}) \phi_n(\mathbf{y}) d\mathbf{y}.
\label{soleig}
\end{eqnarray}
The series expansion allows us to prove the `smoothing effect' for $t>0$: 
$u(t) \in H^2(\Omega)$. In fact, $u(t) \in H^k(\Omega),$ for all $k$.\\

Since $\Omega$ is a bounded domain, a $L^2_{\mathbf x}$ estimate 
implies a $L^q_{\mathbf x}$ estimate for $q \in [1, 2]$.  In the same way,
a $L^{\infty}_{\mathbf x}$ estimate implies a  $L^q_{\mathbf x}$ estimate 
for $1 \leq q \leq \infty$. 
The decay of the norms of the solutions of the pure initial value problem
is summarized in the following result.

\begin{theo} {\bf (Decay for the initial value Neumann problem).}
\label{decayu}
If $u_0 \in L^{\infty}(\Omega)$, the solution $u$ of 
equations (\ref{ec:uin})-(\ref{ec:inuin}) satisfies: 
\begin{eqnarray}
\|u(t) \|_{L^{\infty}_{\mathbf x}} &\leq&  \|u_0 \|_{L^{\infty}_{\mathbf x}}, 
\label{estinfin} \\
\|u(t) \|_{L^{2}_{\mathbf x}} &\leq&  \|u_0 \|_{L^{2}_{\mathbf x}},
\label{est2in} \\
\|\nabla_{\mathbf x} u(t) \|_{L^{2}_{\mathbf x}} &\leq& 
{1\over t^{1/2}} \|u_0 \|_{L^{2}_{\mathbf x}}, 
\label{estd2in} \\
\|\nabla_{\mathbf x} u(t) \|_{L^{r}_{\mathbf x}} &\leq& 
{C_{rq}\over t^{1/2+N/2(1/q-1/r)}} 
\|u_0 \|_{L^{q}_{\mathbf x}}, \quad 1 \leq q \leq r \leq \infty, C_{rq}>0, 
\label{estdqpin} 
\end{eqnarray}
for $t\in (0,T]$, $T>0$.
Moreover, if $\nabla_{\mathbf x} u_0 \in L^{2}_{\mathbf x}$
and $\Delta_{\mathbf x} u_0 \in L^{\infty}_{\mathbf x}$, then
\begin{eqnarray}
\|\nabla_{\mathbf x} u(t) \|_{L^2_{\mathbf x}} &\leq&  
\|\nabla_{\mathbf x} u_0 \|_{L^2_{\mathbf x}},
\label{estd2din} \\
\|\nabla_{\mathbf x} u(t) \|_{L^{\infty}_{\mathbf x}} &\leq&  
C(\| u_0 \|_{L^{\infty}_{\mathbf x}},
\|\Delta_{\mathbf x} u_0 \|_{L^{\infty}_{\mathbf x}}). 
\label{estdinfin} 
\end{eqnarray}
\end{theo}

{\bf Proof.}

By Proposition \ref{comparisonneuman}, $u$ is bounded from above and below 
by the solutions of equations (\ref{ec:uin})-(\ref{ec:inuin}) with initial data 
$\|u_0\|_{\infty}$ and $-\|u_0\|_{\infty},$ respectively. This proves (\ref{estinfin}). 

Multiplying equation (\ref{ec:uin}) by $u$ and integrating, we find the energy 
identity:
\begin{eqnarray}
{1\over 2} \|u(t)\|_2^2 + d \!\! \int_0^t \!\!\int_{\Omega} |\nabla_{\mathbf x} u|^2 =   
{1\over 2} \|u_0\|_2^2,   \label{energy}
\end{eqnarray}
which implies estimate (\ref{est2in}). 

To prove inequality (\ref{estd2in}) we argue by density, assuming first that $u_0\in H^1(\Omega)$. We multiply equation (\ref{ec:uin}) by $u_t$ and integrate 
over $\Omega$ to get
\begin{eqnarray*}
 \|u_t(t)\|_2^2 + {d\over 2} {d \over dt} \!\!\int_{\Omega} |\nabla_{\mathbf x} u(t)|^2 = 0.   
\end{eqnarray*}
We conclude that $\!\!\int_{\Omega} |\nabla u(t)|^2$ 
decreases with time. Inserting this information in identity (\ref{energy}), we find:
\begin{eqnarray*}
{1\over 2} \|u(t)\|_2^2 + t \int_{\Omega} |\nabla_{\mathbf x} u(t)|^2 \leq  
{1\over 2} \|u_0\|_2^2 \; \Rightarrow \;
\|\nabla_{\mathbf x} u(t) \|_{L^{2}_{\mathbf x}} 
\leq {1\over t^{1/2}} \|u_0 \|_{L^{2}_{\mathbf x}}. 
\end{eqnarray*}
The inequality extends to $u_0 \in L^2_{\mathbf x}$ by density. 

% nabla u   (ur, utheta/r) en base polar (er, etheta), er=(x/r,y/r) etheta=(-y/r,x/r)
%  (ur, utheta/r, uphi/r/sinthetat  
%  er funcion de sin cos theta, sin cos phi

To prove the $L^r-L^q$ estimate on the gradients, we resort to expressions of the solutions in terms of heat kernels \cite{heatkernel} and  pointwise estimates of the  kernels \cite{chino}. Our annular domain $\Omega$ is not convex, therefore we can only apply results valid for $C^2$ compact manifolds.
In terms of the heat kernel for the Neumann problem, the solution of 
equations (\ref{ec:uin})-(\ref{ec:inuin}) reads:
\begin{eqnarray*}
u({\mathbf x},t)= \int_{\Omega} K({\mathbf x}, {\mathbf y}, t) u_0({\mathbf y}) 
d{\mathbf y},
\quad K({\mathbf x}, {\mathbf y}, t) = \sum_{n\geq 1} e^{-\lambda_n t}
\phi_n({\mathbf x})\phi_n({\mathbf y}),
\end{eqnarray*}
where $\phi_n$ and $\lambda_n$ are the eigenvalues and eigenfunctions of  the homogeneous Neumann problem, see reference \cite{heatkernel}, pp. 104-106. 
 The kernel function $K$ is positive, symmetric
 in the ${\mathbf x}$ and ${\mathbf y}$ variables, and satisfies
 $\int_{\Omega} K({\mathbf x}, {\mathbf y}, t) d{\mathbf y} =1.$ It is the solution
 of a Neumann problem with measure valued initial data 
 $\delta_{\mathbf x}({\mathbf y}).$ For compact Riemannian manifolds
 with $C^2$ smooth boundary, the gradient of the heat kernel satisfies:
\begin{eqnarray*}
|\nabla K({\mathbf x}, {\mathbf y}, t)| \leq {C \over t^{(N+1) \over 2}}
e^{- {\rho({\mathbf x}, {\mathbf y})^2 \over ct }}, \quad t>0, {\mathbf x}, 
{\mathbf y} \in \Omega,
\end{eqnarray*}
for some positive constants $C,c$, where $N$ is the dimension, and $\rho$ the
Riemannian distance. Our domain $\Omega$ is a ring in ${\mathbb R}^N$.
We may find a constant $d'$ such that $\rho({\mathbf x}, {\mathbf y}) \geq 
d' |{\mathbf x} - {\mathbf y}|$.  
Extending the upper bound to an integral over the whole space:
\begin{eqnarray*}
|\nabla u({\mathbf x},t)|= |\int_{\Omega} \nabla K({\mathbf x}, {\mathbf y}, t) 
u_0({\mathbf y}) d{\mathbf y}|
\leq {C \over t^{(N+1) \over 2}} \int_{\mathbb{R}^N}
e^{- {d' |{\mathbf x},- {\mathbf y}|^2 \over ct }} |u_0({\mathbf y})| d{\mathbf y},
\end{eqnarray*}
the $L^r$-$L^q$ estimates (\ref{estdqpin}) follow from standard $L^r$-$L^q$
estimates for solutions of the heat equation in the whole space \cite{giga}.

Expressing the solution in terms of eigenfunctions (\ref{soleig}), the 
$L^2_{\mathbf x}$ norm of the gradient becomes:
\begin{eqnarray}
\int_{\Omega} |\nabla u(\mathbf{x},t)|^2 d \mathbf{x}  
= \sum_{n,m\geq 1} u_{0,n} u_{0,m}  e^{-\lambda_n t} e^{-\lambda_m t}
\left[ \int_{\Omega}  \nabla \phi_n(\mathbf{x}) \nabla \phi_m(\mathbf{x})  
 d \mathbf{x} \right] \nonumber \\
= 
\sum_{n,m\geq 1} u_{0,n} u_{0,m}  e^{-\lambda_n t} e^{-\lambda_m t}
{\lambda_m\over d} \left[ \int_{\Omega}  \phi_n(\mathbf{x})  \phi_m(\mathbf{x})  
 d \mathbf{x} \right] 
\nonumber \\
= 
\sum_{n\geq 1} u_{0,n}^2   e^{-2\lambda_n t}   \lambda_m
\leq 
\sum_{n\geq 1} u_{0,n}^2    {\lambda_m \over d}
= \int_{\Omega} |\nabla u_0(\mathbf{x})|^2 d \mathbf{x},
\end{eqnarray}
after integrating by parts, using definition (\ref{eigenvalue}) and the orthogonality of the eigenfunctions. This proves estimate (\ref{estd2din}).

To estimate the $L^\infty_{\mathbf x}$ norm of the gradient, we notice that differentiating  formula (\ref{soleig}) and assuming ${\partial u_0
\over \partial {\mathbf n}}$ we find:
\begin{eqnarray*}
\Delta u(\mathbf{x},t) 
= \sum_{n\geq 1} u_{0,n} \Delta \phi_n(\mathbf{x}) e^{-\lambda_n t}
= - \sum_{n\geq 1} u_{0,n} {\lambda_n \over d} \phi_n(\mathbf{x}) e^{-\lambda_n t}=
\nonumber \\
\sum_{n\geq 1}  
\left[\int_{\Omega} \hskip -2mm 
u_0(\mathbf{y}) \Delta \phi_n(\mathbf{y}) d\mathbf{y} \right]\!
\phi_n(\mathbf{x}) e^{-\lambda_n t} \!=\!
\sum_{n\geq 1}  
\left[\int_{\Omega} \hskip -2mm
\Delta u_0(\mathbf{y}) \phi_n(\mathbf{y}) d\mathbf{y} \right]\!
\phi_n(\mathbf{x}) e^{-\lambda_n t}.
\end{eqnarray*}
This expression defines a solution of
\begin{eqnarray*}
\tilde u_t(\mathbf x,t)= d \Delta_{\mathbf x} \tilde u(\mathbf x,t),  
&  {\mathbf x} \in \Omega, \, t>0, \\
 {\partial \tilde u \over \partial {\bf n}}(\mathbf x,t) =  0, & \,
\mathbf x \in S_{r_0}\cup S_{r_1}, \,  t>0,  \\
\tilde u(\mathbf{x},0) =\Delta u_0(\mathbf{x}),  & \quad  {\mathbf x} \in \Omega.
\end{eqnarray*}
The comparison principle in Proposition \ref{comparisonneuman} 
yields $\|Ê\Delta_{\mathbf x} u\|_{L^\infty_{\mathbf x}} \leq 
\|Ê\Delta u_0\|_{L^\infty_{\mathbf x}}.$
This inequality extends to $u_0 \in W^{2,\infty}(\Omega)$ by
density.
Finally, Gagliardo-Nirenberg's inequalities \cite{brezis} provide an 
estimate on the gradients:
$\|Ê\nabla_{\mathbf x} u  \|_{L^\infty_{\mathbf x}}  \leq
C(\|Ê\Delta_{\mathbf x} u\|_{L^\infty_{\mathbf x}}, \|u\|_{L^\infty_{\mathbf x}})
\leq C(\|Ê\Delta_{\mathbf x} u_0\|_{L^\infty_{\mathbf x}}, \|u_0\|_{L^\infty_{\mathbf x}}).$
This proves inequality (\ref{estdinfin}) and concludes the proof.\\

Let us now consider the diffusion problem with a source but zero initial
and boundary values:
\begin{eqnarray}
u_t(\mathbf x,t)= d \Delta_{\mathbf x} u(\mathbf x,t) + h(\mathbf x,t) ,  & \quad {\mathbf x} \in \Omega, t>0, \label{ec:uso} \\
{\partial u \over \partial {\bf n}}(\mathbf x,t)= 0,  & \mathbf x \in S_{r_0}\cup S_{r_1},  t>0, \label{ec:bcuso} \\
u(\mathbf{x},0) =0,  & \quad  {\mathbf x} \in \Omega.
\label{ec:inuso}
\end{eqnarray}
For any $h \!\in\! L^{\infty}(0,T;L^2(\Omega))$, there is a unique global solution 
$u \!\in\! C([0,T],L^2(\Omega)),$ see reference
\cite{galerkin}. It is given by the series expansion:
\begin{eqnarray}
u(\mathbf{x},t) =  \sum_{n\geq 1} \phi_n(\mathbf{x}) \int_0^t h_n(s) 
e^{-\lambda_n (t-s)} ds, 
\quad
h_{n}(s)  = \int_{\Omega}  h({\mathbf y},s) \phi_n(\mathbf{y}) d\mathbf{y}.
\nonumber 
\end{eqnarray}
The series expansion implies again the `smoothing effect': $u(t) \in H^2(\Omega)$ 
for $t>0$. In fact, $u(t)\in H^k(\Omega),$ for all $k$ and $t>0$. This solution 
can be rewritten using the semigroup formalism \cite{pazy}. The initial value problem (\ref{ec:uin})-(\ref{ec:inuin}) defines a semigroup $S(t)u_0= u(t)$, 
$u$ being the solution of equations (\ref{ec:uin})-(\ref{ec:inuin}). The solution 
of an inhomogeneous initial value problem with initial datum $u_0$ and source $h$ can be expressed as:
\begin{eqnarray}
u(t)= S(t)u_0 + \int_0^t S(t-s)h(s) ds. \label{semigroup}
\end{eqnarray}
Theorem \ref{decayu} establishes decay estimates for the semigroup
$S(t)$ and its derivatives, applied to different types of initial data. We
can exploit those estimates to infer the decay of the integral term 
representing solutions with a source. The following estimates hold:

\begin{prop}{\bf (Decay for the inhomogeneous problem).}
\label{decayuso}
For sources $h \in L^{\infty}([0,T] \times \Omega)$, the solution $u$ of 
equations (\ref{ec:uso})-(\ref{ec:inuso}) satisfies: 
\begin{eqnarray}
\|u(t) \|_{L^{\infty}_{\mathbf x}}  \leq   t \|h \|_{L^{\infty}_t L^{\infty}_{\mathbf x}}, 
\qquad
\|\nabla_{\mathbf x} u(t) \|_{L^{\infty}_{\mathbf x}} \leq   2t^{1/2}
\| h \|_{L^{\infty}_t L^{\infty}_{\mathbf x}}, \\
\|u(t) \|_{L^{2}_{\mathbf x}}  \leq   t \|h \|_{L^{\infty}_t L^{2}_{\mathbf x}},
\qquad \;\;
\|\nabla_{\mathbf x} u(t) \|_{L^{2}_{\mathbf x}}  \leq   
 2 t^{1/2}   \|h \|_{L^{\infty}_t L^{2}_{\mathbf x}}, \;\;
\label{estd2so} \\
\|\nabla_{\mathbf x} u(t) \|_{L^{r}_{\mathbf x}}  \leq  C_{rq}  t^{1/2-N/2(1/q-1/r)}  
\|h \|_{L^{\infty}_t L^{q}_{\mathbf x}}, \quad  C_{rq}>0,  {1\over N} > {1\over q}
- {1\over r} >0,
\label{estdqpso} 
\end{eqnarray}
for $t \in [0,T]$.
\end{prop}

{\bf Proof.} Consequence of the semigroup expression (\ref{semigroup}) for the
solutions and Proposition \ref{decayu}. \\

Let us now apply the previous decay estimates to solutions $c$ of  equations (\ref{eq:c})-(\ref{eq:pc0}),(\ref{bc:c}). Let $c_b$ be a function such that $c_b=c_{r_0}$ on $r=r_0$ and $c_b=0$ on $r=r_1$. For simple choices of $c_{r_0}$ this can be done explicitly. Otherwise, we may resort to solutions $c_b$ of the boundary value problem for the
heat equation, with zero initial data, zero source
term, and non homogeneous boundary conditions $c_b=c_{r_0}$ on $r=r_0$ and $c_b=0$ on $r=r_1$. Existence of such solutions has been established in reference \cite{brown} by the method of layer potentials  for boundary data satisfying integrability conditions that always hold for bounded data.
We set $ c = \tilde{c} + c_b$. Then, 
\begin{eqnarray}
\tilde{c}_t- d \Delta_{\mathbf x} \tilde{c} = - \eta c j - c_{b,t} 
+ d\Delta_{\mathbf x} c_{b},  & \quad {\mathbf x} \in \Omega, t>0, \label{ec:ctilde} \\
{\partial \tilde{c} \over \partial {\bf n}}(\mathbf x,t)= 0,  
& \mathbf x \in S_{r_0}\cup S_{r_1},  t>0,  \label{ec:bctilde} \\
\tilde{c}(\mathbf{x},0) =c_0(\mathbf{x})-c_b(\mathbf{x},0),  & \quad  {\mathbf x} \in \Omega.
\label{ec:intilde}
\end{eqnarray}
The term $z=-c_{b,t}+d \Delta_{\mathbf x} c_{b}$ appearing in the right hand side may vanish when $c_b$ is chosen to be a solution of the heat equation. We have the following estimates.

\begin{theo}
\label{estimatesc}
Let $c$ be a solution of equations (\ref{eq:c})-(\ref{eq:pc0}),(\ref{bc:c}) with initial and boundary data verifying $c_0 \in W^{2,\infty}(\Omega)$, $c_0\geq 0$ and %$c_{r_0}\in L^{\infty}(0,T; W^{2,\infty}(\partial \Omega)),$ 
$c_{r_0}\in L^{\infty}(0,T; L^{\infty}(\partial \Omega)),$ 
$T>0$. Let $c_b \in W^{1,\infty}(0,T;W^{2,\infty}(\Omega))$ be a function satisfying $c_b=c_{r_0}$ on $r=r_0$ and $c_b=0$ on $r=r_1$.
Set $K=\mbox{max}(\|c_0-c_b\|_{\infty},\|c_{b,t} 
-d\Delta c_{b}\|_{\infty}). $ Then, $c\geq 0$ and
\begin{eqnarray}
\|c(t)\|_q  \leq [\|c_b\|_{\infty}+K (1+T)] \, \mbox{meas}(\Omega)^{1/q},\quad
t \in [0,T], 1\leq q \leq \infty. \label{estcq}
\end{eqnarray}
Moreover,
\begin{eqnarray}
&\|\nabla c(t) \|_{\infty} \leq \|\nabla c_b (t) \|_{\infty}
+ C( \| c_0-c_b(0) \|_{\infty},   \| \Delta c_0- \Delta c_b(0) \|_{\infty}) 
\label{estdcinf} \\
&+  2t^{1/2} \| c_{b,t} \!-\!d \Delta c_{b}  \|_{\infty}
+ \eta C_{q} t^{{1\over 2}-{N\over 2q}} \|c j \|_{L^{\infty}_t L^{q}_{\mathbf x}}, 
\quad q>N, \nonumber  
\end{eqnarray}
\begin{eqnarray}
&\|\nabla c(t) \|_{2} \leq \|\nabla c_b(t) \|_{2} 
+ \|\nabla c_0\!-\!\nabla c_b(0)\|_{2} 
\label{estdc2} \\
&+2 t^{1/2}  \| c_{b,t}\!-\!d\Delta c_{b} \|_{L^{\infty}_t L^{2}_{\mathbf x}} 
+2  t^{1/2}   \eta\|c j \|_{L^{\infty}_t L^{2}_{\mathbf x}},
\nonumber  
\end{eqnarray}
\begin{eqnarray}
&d \int_0^t  \int_{\Omega}  |\nabla c(s)|^2 ds d{\mathbf x} 
\leq d \int_0^t \int_{\Omega}  |\nabla c_b(s)|^2 ds d{\mathbf x} 
+ {1\over 2}\|  c_0\!-\!  c_b(0)\|_{2}^2
 \label{estdc22} \\
&+ \|c_{b,t}-d\Delta c_{b}\|_{L^2(0,t;L^2_{\mathbf x})} 
\|c\|_{L^2(0,t;L^2_{\mathbf x})}.
\nonumber 
\end{eqnarray}
\end{theo}

{\bf Proof.}
By the comparison principle \ref{comparisonneuman}, we know that $c \geq 0$ and that $ \tilde{c} $ is bounded from above by the solution $\tilde{C}$ of system (\ref{ec:ctilde})-(\ref{ec:intilde}) with right hand side $z=-c_{b,t} +d\Delta c_{b}$. By Proposition \ref{comparisonneuman}, $|\tilde{C}|\leq K(1+t).$ Since $\Omega$ is a bounded domain, estimate (\ref{estcq}) follows.

Let us now study the derivatives of $\tilde{c}$. The energy inequality provides
a uniform $L^2_{\mathbf x t}$ estimate that implies inequality (\ref{estdc22}):
\begin{eqnarray*}
{1\over 2}\|\tilde{c}(t)\|_{L^2_{\mathbf x}}^2 + d \int_0^t\int_{\Omega} |\nabla \tilde{c}(s)|^2 ds d{\mathbf x} \leq  {1\over 2}\|\tilde{c}(0)\|_{L^2_{\mathbf x}}^2
\nonumber \\
+ \|d\Delta c_{b}-c_{b,t}\|_{L^2(0,T;L^2_{\mathbf x})} \|c\|_{L^2(0,T;L^2_{\mathbf x})}.
\end{eqnarray*}
%
% - \int \Delta c c = \int |\nabla c|^2 - \int d_n c c
%
To prove inequalities (\ref{estdc2}) and (\ref{estdcinf}), we split $\tilde c = \tilde c_1
+ \tilde c_2$. By linearity,  we take $\tilde{c}_1$ to be a solution of a heat equation with initial datum $\tilde{c}_0$, zero source and zero boundary condition. We choose $\tilde{c}_2$ to be the solution of another heat problem, with source $- \eta c j + z$, plus zero initial and boundary conditions. 
The estimates stated in Theorem \ref{decayu} hold for $\tilde{c}_1$ and those in Proposition \ref{decayuso} to $\tilde{c}_2$. Differentiating $- \eta c j $ must be avoided, not to introduce the spatial derivatives of $c$ we intend to control. In this way, we obtain inequalities (\ref{estdc2}) and (\ref{estdcinf}). 

\section{Nonlinear problem with known boundary condition}
\label{sec:nonlinearbc}

Solutions for the the nonlinear angiogenesis model may be constructed
employing an iterative scheme. For $m\geq 2$, we consider the 
linearized system of equations
\begin{eqnarray} 
\frac{\partial}{\partial t} p_m(\mathbf{x},\mathbf{v},t) \!+\! 
\mathbf{v} \!\cdot\! \nabla_\mathbf{x}   p_m(\mathbf{x},\mathbf{v},t) 
\!+\!  \nabla_\mathbf{v} \!\cdot\! [({\mathbf F}(c_{m-1}(\mathbf{x},t)) 
\!-\! \beta \mathbf{v} ) p_m(\mathbf{x},\mathbf{v},t) ]  
 \label{eq:pm} \\
\!-\! \sigma \Delta_\mathbf{v} p_m(\mathbf{x},\mathbf{v},t)
\!+\! \gamma b_{m-1}(\mathbf{x},t)p_{m}(\mathbf{x},\mathbf{v},t) 
\!=\! \alpha(c_{m-1}(\mathbf{x},t)) \nu(\mathbf{v}) p_{m}(\mathbf{x},\mathbf{v},t),  \nonumber  \\
b_{m-1}(\mathbf{x},t) = \int_0^t \! \!ds \! \! \int_{\mathbb R^N} \! \! d{\bf v}' 
p_{m-1}(\mathbf{x},\mathbf{v}',s),   \label{eq:am} \\
 \alpha(c_{m-1})=\alpha_1\frac{c_{m-1}}{c_R+ c_{m-1}}, \quad
 {\bf F}(c_{m-1})= \frac{d_1}{(1+\gamma_1c_{m-1})^{q_1}}
 \nabla_{\mathbf x} c_{m-1},
\label{eq:alphaFm} \\
p_m(\mathbf{x},\mathbf{v},0) = p_0(\mathbf{x},\mathbf{v}),  \label{eq:pm0} \\
\frac{\partial}{\partial t}c_{m-1}(\mathbf{x},t) = d \Delta_{\mathbf x} c_{m-1}(\mathbf{x},t) 
- \eta c_{m-1}(\mathbf{x},t) j_{m-1}(\mathbf{x},t), \label{eq:cm} \\
j_{m-1}(\mathbf{x},t) \!=\! \int_{\mathbb R^N}  {|\mathbf{v}| \over 1 + e^{|\mathbf{v}- \mathbf{v}_0 \chi |^2 / \sigma_v^2}}
p_{m-1}(\mathbf{x},\mathbf{v},t)\, d \mathbf{v}  \label{eq:jm}, \\
c_{m-1}(\mathbf{x}) = c_0(\mathbf{x},0),   \label{eq:cm0}
\end{eqnarray}
supplemented with the boundary conditions:
\begin{eqnarray}
& {\partial c_{m-1} \over \partial {\bf n}}({\mathbf x},t) = c_{r_0}({\mathbf x},t)<0, 
\; {\mathbf x} \in S_{r_0},  \quad
 {\partial c_{m-1} \over \partial {\bf n}} ({\mathbf x},t) = 0, \;
{\mathbf x} \in S_{r_1}, \label{bc:cm} \\
& p_m({\mathbf x},{\mathbf v},t) = g({\mathbf x},{\mathbf v},t)  \geq 0,  \quad
\mbox{for} \; \mathbf v \cdot \hat {\mathbf n}<0, \; {\mathbf v} \in {\mathbb R}^N, \; 
{\mathbf x} \in S_{r_0}\cup S_{r_1}. \label{bc:rm}  
\end{eqnarray}

We initialize the scheme setting $p_1=0$ and ${j}_1=0$. 
$c_1$ is the solution of the associated heat equation.  The function
$p_2$ is a nonnegative solution of the Fokker-Planck problem with smooth 
and bounded coefficient fields ${\mathbf F(c_1)}$ and $\alpha(c_1)$
in a bounded domain. Let us see that 
the resulting sequence is well defined under our hypotheses on the
data and we may extract a subsequence converging
to a solution of the original problem.

\begin{theo} 
\label{maintheorem}
Let us assume that
\begin{eqnarray}
p_0 \geq 0, c_0 \geq 0, g \geq 0, \label{signo} \\
c_0  \in W^{2,\infty}(\Omega),  \label{c0inf} \\
(1+|{\mathbf v}|^2)^{\mu/2} p_0
\in L^{\infty}\cap L^1(\Omega \times {\mathbb R}^N), \quad \mu >N,
\label{p0inf} \\
c_{r_0}\in L^{\infty}(0,T; L^{\infty}(S_{r_0})), \label{bcm} \\
(1+| {\mathbf v} \cdot {\mathbf n}|) (1+|{\mathbf v}|^2)^{\mu/2} g
\in {L^{\infty}(0,T;L^\infty\cap L^1(\Gamma^-))}, \label{bpm}
\end{eqnarray}
and that a function $c_b$ is found verifying the hypotheses of Theorem \ref{estimatesc}.
Then, there exists a nonnegative solution $(p, c)$ of equations 
(\ref{eq:c})-(\ref{bc:plin}) satisfying:
\begin{eqnarray}
c \in L^{\infty}(0,T;W^{1,\infty}(\Omega)), \label{dcinf} \\
p \in L^{\infty}(0,T;L^{\infty}\cap L^1(\Omega \times {\mathbb R}^N)),
\nabla_{\mathbf v }p \in L^{2}(0,T;L^2(\Omega \times {\mathbb R}^N)),
\label{pinf} \\
(1+|{\mathbf v}|^2)^{\mu/2}p \in 
L^{\infty}(0,T;L^{\infty}\cap L^1(\Omega \times {\mathbb R}^N)),  
\label{pdecay} \\
p \in L^{\infty}(0,T;L^\infty_{\mathbf x}(\Omega,
L^1_{\mathbf v}( {\mathbb R}^N)), \label{pmominf}  
\end{eqnarray}
for any $T>0$. 
\end{theo}

The proof is organized in several steps. First, we argue that the scheme is well defined. Then, we obtain  uniform estimates on the $L^q$ norms of the solutions of the iterative scheme. Next, we derive $L^{\infty}$ estimates on the coefficients ${\mathbf F}_{m-1}$, $j_{m-1}$ and $b_{m-1}$ using the velocity decay. Estimates on the derivatives of the densities with respect to ${\mathbf v}$ allow us to pass to the limit in the equations using compactness results for the specific of FP operator, obtaining a nonnegative solution of the nonlinear problem with the stated regularity. \\

\noindent {\bf Proof.}\\
\leftline{\it Step 1: Existence of nonnegative solutions for the scheme.}
First, let us argue that the scheme (\ref{eq:pm})--(\ref{bc:rm}) is well defined.
Setting $p_1=0$, we have ${j}_1=0$ and due to (\ref{eq:cm}),
$c_1(\mathbf{x},t)$ is the solution of the associated heat equation
\begin{eqnarray*}
\frac{\partial}{\partial t}c_{1}(\mathbf{x},t) = d \Delta_{\mathbf x} c_{1}(\mathbf{x},t)
\quad \mathbf x \in \Omega, t>0, \\ 
c_{1}(\mathbf{x},0) = c_0(\mathbf{x}), \quad \mathbf{x} \in \Omega, \\
{\partial c_{1} \over \partial {\bf n}}({\mathbf x},t) = c_{r_0}({\mathbf x},t)<0, 
\; {\mathbf x} \in S_{r_0},  \quad
 {\partial c_{1} \over \partial {\bf n}} ({\mathbf x},t) = 0, \;
{\mathbf x} \in S_{r_1}, \quad t>0,
\end{eqnarray*}
satisfying the properties of Theorem \ref{estimatesc}. The function $p_2$ 
is the nonnegative solution of the Fokker-Planck problem with smooth
and bounded coefficient fields ${\mathbf F(c_1)}$ and $\alpha(c_1)$
in a bounded domain, i.e.,
$$
\alpha(c_{1})=\alpha_1\frac{c_{1}}{c_R+ c_{1}}, \quad
{\bf F}(c_{1})= \frac{d_1}{(1+\gamma_1c_{1})^{q_1}}\nabla_{\mathbf x} c_{1}.
$$

Let us proceed by induction. We assume that ${j}(p_{m-1})$ 
and $b(p_{m-1})$ are nonnegative bounded functions. Then, $c_{m-1}$ is the 
unique solution of equations (\ref{eq:cm})-(\ref{eq:cm0}) with boundary conditions (\ref{bc:cm}), whose existence can be proven by Galerkin or spectral methods \cite{galerkin}. By Proposition \ref{comparisonneuman} we know that 
$c_{m-1} \geq 0$ if $c_0\geq 0$.  This implies that 
$0\leq \alpha(c_{m-1}) \leq \alpha_1$ and
$0 \leq \frac{d_1}{(1+\gamma_1c_{m-1})^{q_1}} \leq d_1.$ 
Moreover,  Theorem \ref{estimatesc} provides $L^\infty$ and $L^q$ 
bounds for $c_{m-1}$.
Then, Theorem \ref{estimatesc} implies that 
$\nabla_{\mathbf x} c_{m-1}$ is a bounded function and also
${\mathbf F}(c_{m-1})$.
Since $\alpha(c_{m-1})$ and ${\mathbf F}(c_{m-1})$ are bounded, and
$b(p_{m-1})$ is assumed to be bounded, $p_m$ is the unique nonnegative
solution of equations (\ref{eq:pm}),(\ref{eq:pm0}) with boundary conditions
(\ref{bc:rm}) that satisfies the estimates collected in 
Theorems \ref{positivity}, \ref{conservation} and \ref{comparison}.
This implies that $\gamma b(p_{m-1}) p_{m} \geq 0$ and 
$\alpha(c_{m-1}) \nu \leq  \alpha_1\| \nu \|_{\infty}$. 
By Lemma \ref{cotaj}, $\|j(p_m)\|_{L^\infty_{\mathbf x}}$ is bounded 
in terms of $\|p_m\|_{L^\infty_t L^\infty_{\mathbf x \mathbf v}}.$
Proposition \ref{velocityLinf} implies that $b(p_m) \in 
L^{\infty}(\Omega \times (0,T))$. This allows us to construct $c_m$
and $p_{m+1}$, and so on. \\

\leftline{\it Step 2: A priori estimates on the tumor angiogenic factor $c_m$.}

By Theorem \ref{estimatesc}, setting 
$K=\mbox{max}(\|c_0-c_b\|_{L^{\infty}_{\mathbf x}},\|c_{b,t} -d\Delta c_{b}\|_{L^{\infty}_{\mathbf x t}}),$ we get
\begin{eqnarray}
\|c_m(t)\|_{L^{q}_{\mathbf x}}  \leq (\|c_b\|_{\infty}+K T) 
\, {\rm meas}(\Omega)^{1/q},\quad
t \in [0,T], 1\leq q \leq \infty. \label{uniformcinf}
\end{eqnarray}

% int Delta c c^{q-1} = - int nabla c (q-1) c^{q-2} nabla c
%  + int dc/dn c^{q-1}
The energy inequality yields a bound on the gradient independent
of ${j}(p_{m-1})$:
\begin{eqnarray*}
{d\over dt} \int_{\Omega} c^2_m d{\mathbf x} + 
d \int_{\Omega}   |\nabla_{\mathbf x} c_m|^2  d{\mathbf x} 
\leq d \int_{\partial \Omega} c_{r_0} c_m.
\end{eqnarray*}
Integrating in time we find
\begin{eqnarray}
d \int_0^T \int_{\Omega}   |\nabla_{\mathbf x} c_m|^2  d{\mathbf x} 
\leq \|Êc_0\|_{L^2_{\mathbf x}}^2 
+ d \|c_{r_0}\|_{L^1(\partial \Omega \times (0,T))} \|c_m\|_{\infty}.
\label{uniformdc2}
\end{eqnarray}
Theorem \ref{estimatesc} provides  alternative energy estimates (\ref{estdc22}) on the $L^2_t L^2_{\mathbf x}$ norm of $\nabla_{\mathbf x} c_{m-1}$ and the $L^{\infty}_{\mathbf x t}$ norm of $\nabla_{\mathbf x} c_{m-1}$:
\begin{eqnarray}
\|\nabla c_{m-1}(t) \|_{\infty} \leq \|\nabla_{\mathbf x} c_b (t) \|_{\infty}
+ C( \| c_0-c_b(0) \|_{\infty},   \| \Delta_{\mathbf x} c_0- \Delta_{\mathbf x} c_b(0) \|_{\infty}) \nonumber \\
+  t^{1/2} \| c_{b,t} \!-\!d \Delta_{\mathbf x} c_{b}  \|_{\infty}
+ \eta C_{q} t^{{1\over 2}-{N\over 2q}} \|c_{m-1}  {j}(p_{m-1})  \|_{L^{\infty}_t L^{q}_{\mathbf x}}, \quad q>N, \label{uniformdcinf} 
\end{eqnarray}
for $t\in [0,T].$ \\

\leftline{\it Step 3: A priori estimates on the tip vessel density $p_m$.}

Let us revisit the $L^q$ estimates in Theorem \ref{conservation}. The
conservation of mass implies inequality (\ref{lp1}) with $q=1$:
\begin{eqnarray}
\|p_m(t)\|_{L^1_{\mathbf x \mathbf v}} \leq
\Big( \|p_0\|_{L^1_{\mathbf x \mathbf v}}  + \int_{\Sigma^-_T} |{\mathbf v} \cdot {\mathbf n}({\mathbf x}) | g \Big) e^{\alpha_1 \|\nuÊ\|_{\infty} t}, \quad t \in [0,T].
\label{uniformL1xv}
\end{eqnarray}
Applying inequality (\ref{lp1}) with $1< q<\infty$, we find:
\begin{eqnarray}
\|p_m(t)\|_{L^q_{\mathbf x \mathbf v}}^q \leq
\Big( \|p_0\|_{L^q_{\mathbf x \mathbf v}}^q  \!\!+\!\! \int_{\Sigma^-_T} |{\mathbf v} 
\cdot {\mathbf n}({\mathbf x}) | g^q \Big) e^{\big( N \beta (q-1) + \alpha_1 \|\nuÊ\|_{\infty} \big) t}, \quad t \in [0,T].
\label{uniformLqxv}
\end{eqnarray}
Uniform $L^{\infty}$ estimates follow either inequality (\ref{lp1}) with $q=\infty$: 
\begin{eqnarray}
\|p_m\|_{L^{\infty}(Q_T)} \leq\Big( \|p_0\|_{L^\infty(\Omega \times \mathbb{R}^N)} 
+ \|g\|_{L^{\infty}_k(\Sigma^-_T)}  \Big)  e^{\big( N \beta + \alpha_1 \|\nuÊ\|_{\infty} \big) t},
\label{uniformLinfxv}
\end{eqnarray}
or Theorem \ref{comparison} with $\|g\|_{L^{\infty}_k(\Sigma^-_T)}$ replaced
by $\|g\|_{L^{\infty}(\Sigma^-_T)}$.

To be able to extract a converging subsequence from the sequence $p_m$, we need estimates on its derivatives. Let us revisit the $L^2$ estimate (\ref{lp}) provided by Theorem \ref{conservation} for $p_m$:
\begin{eqnarray} \hskip -8mm
\frac{d}{d t}\|p_m(t)\|^2_{L^2_{\mathbf x \mathbf v}} \hskip -3mm &=&
\hskip -3mm
\int_{\Gamma_{-}}|{\mathbf v}\cdot {\mathbf n}({\mathbf x})|g^2dS \, d{\mathbf v}
+ \int_{\Omega\times\mathbb{R}^N} \alpha(c_{m-1}) \nu p_m^2 d{\mathbf x} \, 
d{\mathbf v}    \nonumber \\
&-&  \hskip -3mm \int_{\Gamma_{+}}|\mathbf v\cdot \mathbf n(x)|
{(\rm Tr} p_m)^2 dS \, 
d{\mathbf v} -\int_{\Omega\times\mathbb{R}^N} \gamma b(p_{m-1}) p_m^2 
d{\mathbf x} \, d{\mathbf v}   \nonumber \\
&+& \hskip -3mm N\beta \|p_m(t)\|^2_{L^2_{\mathbf x \mathbf v}} 
\hskip -2mm - 2 \sigma  \!\!\int_{\Omega\times \mathbb{R}^N} \hskip -3mm 
|\nabla_{\mathbf v} p_m|^2  d{\mathbf x} \, d{\mathbf v}.
\nonumber %\label{energydiff}
\end{eqnarray}
Integrating in time and neglecting negative terms, we find
\begin{eqnarray}  
 2 \sigma    \int_{Q_T}   
|\nabla_{\mathbf v} p_m|^2  d{\mathbf x} \, d{\mathbf v} ds \leq  
\|p_0\|^2_{L^2_{\mathbf x \mathbf v}} \!\! + 
\!\!\int_{\Sigma_{-}^T} \!\!  
|{\mathbf v}\cdot {\mathbf n}({\mathbf x})|g^2 dS d{\mathbf v} ds \nonumber \\
+ (\alpha_1 \|\nu\|_{\infty} + N \beta) \int_{Q_T}  p_m^2 d{\mathbf x} \, d{\mathbf v}   ds.
\nonumber %\label{energyinteg}
\end{eqnarray}
The uniform estimates on $\|p_m\|_{L^2(Q_T)}$ yield
a uniform estimate on $\|Ê\nabla_{\mathbf v} p_m\|_{L^2(Q_T)}$. \\

\leftline{\it Step 4: Uniform bounds on velocity integrals and velocity decay
of $p_m$.}

In Steps 2 and 3 we have obtained uniform estimates on the blood vessel
density norms $\|p_m\|_{L^\infty(0,T;L^{q}_{\mathbf x\mathbf v})}$ 
and the tumor angiogenic factor norms $\|c_m\|_{L^\infty(0,T;L^{q}_{\mathbf x})}$
for $1 \leq q \leq \infty$.

Lemma \ref{cotaj} provides a uniform bound of the $L^q_{\mathbf x}$
norms of the fluxes ${j}_{m-1}$, $1 \leq q \leq \infty$ in terms
of the bounds (\ref{uniformLinfxv}) 
on $\| p_m\||_{L^\infty(0,T;L^q_{\mathbf x\mathbf v})}$ established
in Step 3. Thanks to inequality (\ref{uniformdcinf}) in  Step 2, 
we obtain a uniform estimate on  
$\| \nabla_{\mathbf x} c_{m-1} \|_{L^\infty(0,T;L^\infty(\Omega))}$. 
A uniform estimate on 
$\|\mathbf F (c_{m-1})Ê\|_{L^\infty(0,T;L^\infty(\Omega))}$
follows.

Next, we apply Proposition \ref{velocityLinf} to equation (\ref{eq:pm}), 
setting $a= \gamma b(p_{m-1}) -\alpha(c_{m-1}) \nu$ and
${\mathbf F}={\mathbf F}(c_{m-1})$, with ${j}={j}_{m-1}$ 
depending  on $p_{m-1}$.
Step 1 guarantees that $a \in L^{\infty}$. Its negative part 
$a^-(c_{m-1})=\alpha(c_{m-1})\nu$  satisfies 
$\|a^-\|_{\infty}\leq \alpha_1 \|\nu \|_{L^\infty_{\mathbf v}}.$   
Thanks to the uniform estimate on $\|\mathbf F (c_{m-1})\|_{\infty}$
and $\|p_m\|_{\infty}$, Proposition \ref{velocityLinf} provides
a uniform estimate on 
$\|(1\!+\! |{\mathbf  v}|^2)^{\mu/2}  p_{m-1} \|_{L^\infty_{\mathbf x \mathbf v}}$.
Then, inequality (\ref{boundinf}) in Lemma \ref{interp}  yields a uniform estimate 
on $\|p_m\|_{L^\infty(0,T;L^{\infty}_{\mathbf x} L^1_{\mathbf v})} $.
We also obtain as a consequence an upper bound of the form
\begin{eqnarray}
|p_m| \leq {C \over (1 + |{\mathbf v}|^2)^{\mu/2}} = {\cal P}, \quad \mu > N, C>0.
\label{decayvpm}
\end{eqnarray}
This upper bound ${\cal P}$ is integrable, and belongs to $L^q(\Omega
\times \mathbb R^N)$ for any $q \in[1,\infty]$ since $\mu >N$ and
$\Omega$ is bounded.

In conclusion, the coefficients $b_{m-1}$, $\alpha(c_{m-1})$, ${j}_{m-1}$ 
and $F(c_{m-1})$ appearing in the equations are uniformly bounded in 
$L^{\infty}(0,T; L^{\infty}_{\mathbf x})$. \\

\leftline{\it Step 5: Compactness of the iterates.}

Once we have obtained uniform estimates on $p_m$ and their velocity
derivatives, we resort to the compactness results in reference \cite{bouchut}
to extract converging subsequences. 

\begin{lemma} \cite{bouchut}
 Let $\sigma > 0$, $\beta \geq 0$, $T >0$, $1 \leq q < \infty$,
$p_0 \in L^q({\mathbb R}^N)$, $h \in L^1(0,T;L^q({\mathbb R}^N\times {\mathbb R}^N))$ and
consider the solution $p \in C([0,T], L^q( \mathbb{R}^N \times \mathbb{R}^N))$ 
of:
\begin{eqnarray}
{\partial p \over \partial t} + {\mathbf v} \nabla_{\mathbf x} p
- \beta {\rm div}_{\mathbf v}({\mathbf v}  p) - \sigma \Delta_{\mathbf v} p
&=& h \quad {\rm in} \; \mathbb{R}^N \times \mathbb{R}^N  \times  
(0,T), \label{ecF0} \\
p(0) &=& p_0  \quad {\rm in} \; \mathbb{R}^N \times \mathbb{R}^N. \nonumber 
\end{eqnarray}
Assume that $p_0$ belongs to a bounded subset of $L^q( \mathbb{R}^N \times \mathbb{R}^N)$
and $h$ belongs to a bounded subset of $L^r(0,T;L^q({\mathbb R}^N\times {\mathbb R}^N))$
with $1 < r \leq \infty$. Then, for any $\eta >0$ and any bounded open subset
$\omega$ of ${\mathbb R}^N\times {\mathbb R}^N$, $p$ is compact in $C([\eta,T],L^q(\omega))$.
\end{lemma}
%Previous compactness results by R.J. DiPerna and P.L. Lions \cite{diperna} 
%guarantee compactness in $L^1(0,T;L^q(\omega))$ for more general operators.

These results are stated for problems set in the whole space. Here, we deal
with a problem set in $\Omega \subset {\mathbb R}^N$. We may extend them
to the whole space multiplying by functions $\phi \in C^{\infty}_c(\Omega)$.
The truncated sequences $q_m= \phi p_m$ satisfy
\begin{eqnarray}
{\partial q_{m} \over \partial t} + {\mathbf v} \nabla_{\mathbf x} q_{m}
- \beta {\rm div}_{\mathbf v} ({\mathbf v}  q_{m}) - \sigma \Delta_{\mathbf v} q_{m}  =  h_m \quad {\rm in} \; \mathbb{R}^N \times \mathbb{R}^N  \times  
(0,T),  \nonumber 
\end{eqnarray}
where the sources 
\begin{eqnarray}
h_m= - {\mathbf v}\cdot\nabla_{\mathbf x}  \phi \, p_{m} -
\phi \; {\bf F}(c_{m-1}) \cdot \nabla_{\mathbf v} p_{m} 
-\gamma b(p_{m-1})  p_{m} \phi \,
+ \alpha(c_{m-1}) \nu p_m \phi \nonumber
\end{eqnarray}
are bounded in $L^2(Q_T)$ and the initial state 
$\phi p_0 \in L^1_{\mathbf x \mathbf v} \cap L^{\infty}_{\mathbf x \mathbf v}$ 
is fixed. The sequence $p_m$ is therefore locally
compact and by a diagonal extraction procedure we may extract a 
subsequence $p_{m'}$ converging to a limit $p$ pointwise,
and strongly in $C([\eta,T],L^2(\omega))$ for any 
$\omega \subset \Omega \times \mathbb R^N$. 
Uniform bounds together with uniform control of the velocity
decay allow us to extend compactness up to the borders 
\cite{brezis,unbounded}.
Weak convergences of $p_m$ and $\nabla_{\mathbf v} p_m$ hold
in all the spaces in which we have uniform estimates.

In Step 2, we have obtained a uniform bound on $c_m$
in $L^{2}(0,T;H^1(\Omega))$. Step 4 provides a uniform estimate
on ${j}(p_m)$ in $L^{\infty}(0,T;L^{\infty}(\Omega))$.
Using equation (\ref{eq:cm}), we conclude that
${\partial c_m \over \partial t}$  is bounded in $L^2(0,T;H^{-1}(\Omega))$.
Standard compactness results in reference \cite{lions} yield compactness
for the sequence $c_m$ in $L^2(0,T;L^2(\Omega))$. A subsequence  
$c_{m'}$, converges pointwise and strongly in $L^2$ to a function $c$. 
Weak convergences hold in all the spaces for which uniform bounds 
have been established.\\

\leftline{\it Step 6: Convergence to a solution.}

Let us first pass to the limit in the nonlocal terms using the integrable
upper bound $\cal P$ defined in (\ref{decayvpm}). We know that
$p_{m'}$ and ${\mathbf v} w(\mathbf v)p_{m'}$ converge pointwise to $p$
and $w(\mathbf v) p$. The bounds 
$ 0 \leq p_{m'} \leq {\cal P}\in L^{\infty}(0,T;L^1(\Omega\times {\mathbb R}^N))$ 
and
$|{\mathbf v}| w(\mathbf v)p_{m'} \leq  |\mathbf v| w(\mathbf v){\cal P} \in L^{\infty}(0,T;L^1(\Omega\times {\mathbb R}^N))$
imply pointwise convergence for the nonlocal coefficients:
\begin{eqnarray*}
b(p_{m'}) \rightarrow b(p)\geq 0, \quad {j}(p_{m'}) \rightarrow {j}(p),
\quad a.e. \, {\mathbf x} \in \Omega, t \in [0,T].
\end{eqnarray*}
Let us now consider the nonlinear products. Pointwise convergence
of $b(p_{m'-1})p_{m'}$ to $b(p)p$, together with the bound
$ 0 \leq b(p_{m'-1})p_{m'} \leq b({\cal P}){\cal P}
\in L^1(0,T;L^1(\Omega\times {\mathbb R}^N))$, imply strong convergence
in $L^1(0,T;L^1(\Omega\times {\mathbb R}^N))$. Similarly, pointwise convergence
of $\alpha(c_{m'-1})\nu p_{m'}$ to $\alpha(c) \nu p$, together with the bound
$ |\alpha(c_{m'-1}) \nu p_{m'}| \leq \alpha_1 \|\nu\|_{\infty} {\cal P}
\in L^1(0,T;L^1(\Omega\times {\mathbb R}^N))$, ensure strong convergence
in $L^1(0,T;L^1(\Omega\times {\mathbb R}^N))$.
Finally, pointwise convergence
of $ {j} (p_{m'-1}) c_{m'-1}$ to $ {j} (p) c$, together with the bound
$  {j} (p_{m'-1}) c_{m'-1} \leq {\rm max}_{m'}\| {j} (p_{m'-1})\|_{L^{\infty}(\Omega \times (0,T))} (\|Êc_b\||_{L^{\infty}(\Omega \times (0,T))}+KT) \in L^1(0,T;L^1(\Omega))$, 
yield strong convergence in $L^1(0,T;L^1(\Omega))$.
Strong convergences extend to any $L^q$ with $q$ finite.

The term involving the force field is more complex. Notice that the
sequence $ p_{m'} \frac{d_1}{(1+\gamma_1c_{m'-1})^{q_1}} $ tends  pointwise
to $ p \frac{d_1}{(1+\gamma_1c)^{q_1}} $ and is bounded by ${\cal P} d_1  
\in L^q(0,T;L^q({\Omega\times {\mathbb R}^N}))$ for any $q \in [1,\infty]$. Thus, we have strong convergence in $L^q_{\mathbf x \mathbf v t}$ for all finite $q$. 
The sequence $\nabla_{\mathbf x} c_{m'-1}$ is bounded in
$L^2(0,T;L^2({\Omega\times {\mathbb R}^N}))$. Therefore, it tends
weakly to $\nabla_{\mathbf x} c$ in $L^{2}_{\mathbf x \mathbf v t}$ . 

Using these convergences we pass to the limit in the weak formulation
of the equations for $p_{m'}$:
\begin{eqnarray}  
 \displaystyle
\int_{Q_T} \hskip -3mm p_{m'} \left[\frac{\partial \varphi}{\partial t}
\!+\!  \mathbf v\cdot \nabla_{\mathbf x} \varphi
\!-\! \beta {\mathbf v}\cdot \nabla_{\mathbf v} \varphi
\!+\!  {\bf F}(c_{m'-1}) \cdot \nabla_{\mathbf v} \varphi
\!+\! \sigma\Delta_{\mathbf v} \varphi
\!-\! b(p_{m'-1})\varphi \right] \,d{\mathbf x} d{\mathbf v} dt \nonumber\\ 
\displaystyle
+\int_{\Omega\times \mathbb{R}^N} \hskip -4mm p_0
\varphi({\mathbf x},{\mathbf v},0) \,d{\mathbf x} d{\mathbf v} 
+\int_{\Sigma^-_{T}}|{\mathbf v}\cdot {\mathbf n}({\mathbf x})|
g\varphi \, dS d{\mathbf v} dt
=\int_{Q_T} \hskip -2mm \alpha \nu p_{m'} \varphi \,d{\mathbf x} d{\mathbf v} dt,
  \nonumber \end{eqnarray}
for any $\varphi\in C^\infty(\overline{\Omega} \times \mathbb{R}^N \times [0,T))$ with compact support in ${\mathbf v}$ 
such that $\varphi=0$ on $\Sigma^+_T$.
Weak convergence of $p_{m'}$ is enough to pass to the limit in the linear terms. 
For the rest, we use the strong convergences established above and the
weak convergence of $\nabla_{\mathbf x} c_{m'-1}$. A similar argument can be applied in the weak formulation of equation (\ref{eq:cm}). Therefore, $p$ and $c$ solve the original angiogenesis problem (\ref{eq:c})-(\ref{bc:plin}).

%
% int a_n b_n - int ab = int (an -a) bn + int a (bn-b)
%

%
% Problem for uniqueness: all norm estimates proven for functions
% that are already known to be positive, this is not the case of differences
%

\section{Nonlocal boundary conditions}
\label{sec:nonlocalbc}

In the previous section we constructed solutions for the angiogenesis
model assuming the boundary values for the density known.
The general problem with nonlocal boundary condition 
(\ref{bc:Nr0})-(\ref{bc:Nr1}) becomes  more complex.
Let us address first the linear problem with nonlocal boundary 
conditions. We define the functions:
\begin{eqnarray}
  K_1(\beta,\sigma,\chi,\sigma_v) =  {\rm Max}_{\{  {\mathbf v} \cdot \mathbf n  >0 \} }
 |{\mathbf v} \cdot {\mathbf n}| \; e^{- {\beta \over \sigma}| {\mathbf v} - \mathbf v_0 |^2}  \Big[ \int_{ \{ {\mathbf v} \cdot \mathbf n  < 0 \} } 
 \hskip -4mm {  |\tilde{\mathbf v} \cdot {\mathbf n}| \;
e^{- {\beta \over \sigma }| \tilde{\mathbf v} \!-\!  \mathbf v_0|^2}  \;  
d \tilde{\mathbf v}
\over 1 \!+\! e^{|\tilde{\mathbf v} - \chi  \mathbf v_0|^2/\sigma_v^2}}
\Big]^{-1}, \nonumber
\end{eqnarray}
\begin{eqnarray}
  K_2(\chi,\sigma_v)= \int_{ \{  {\mathbf v} \cdot \mathbf n  >0 \} } {d  {\mathbf v} \over 1 \!+\! e^{| {\mathbf v} - \chi  \mathbf v_0|^2/\sigma_v^2}}.
\label{K1K2}
\end{eqnarray}
For $\chi |\bf v_0| >>1 $ fixed, $  K_2<1$ choosing $\sigma_v$  small enough.
Then, $ K_1 <1$, $ K_1  K_2<1$ choosing 
${\beta \over \alpha} $ small. 
% but \sigma_v cannot go to zero is bounded from below. 

\begin{theo} 
\label{maintheorembc}
Let us assume that
\begin{eqnarray}
p_0 \geq 0, \quad (1+|{\mathbf v}|^2)^{\mu/2} p_0
\in L^{\infty}\cap L^1(\Omega \times {\mathbb R}^N),\quad  \mu >N,
\label{p0inflin}  \\
a \in L^\infty(Q_T), \quad {\mathbf F} \in L^{\infty}(\Omega \times (0,T)),
\label{coefin}  \\
j_0 \geq 0, \quad  j_0 \in L^\infty(\Sigma_T^-).
\end{eqnarray}
Then, there exists a nonnegative solution $p$ of the linear equations 
 (\ref{lin:ec})-(\ref{lin:in}) with boundary conditions (\ref{bc:Nr0})-(\ref{bc:Nr1})  satisfying:
\begin{eqnarray}
p \in L^{\infty}(0,T;L^{\infty}\cap L^1(\Omega \times {\mathbb R}^N)),
\nabla_{\mathbf v }p \in L^{2}(0,T;L^2(\Omega \times {\mathbb R}^N)),
\label{pinflin} \\
(1+|{\mathbf v}|^2)^{\mu/2}p \in 
L^{\infty}(0,T;L^{\infty}\cap L^1(\Omega \times {\mathbb R}^N)),  
\label{pdecaylin} \\
(1+ |\mathbf v \cdot \mathbf n|)(1+|{\mathbf v}|^2)^{\mu/2} {\rm Tr \,}p^\pm \in 
L^{\infty}(0,T;L^{\infty}\cap L^1(\Sigma_T^\pm)),  
\label{ptracelin} \\
p \in L^{\infty}(0,T;L^\infty_{\mathbf x}(\Omega,
L^1_{\mathbf v}( {\mathbb R}^N)), \label{pmominflin}  
\end{eqnarray}
for any $T>0$, provided the parameters $\beta, \sigma, \sigma_v, \chi$  satisfy  $K_1(\beta,\sigma,\chi,\sigma_v)K_2(\chi,\sigma_v)<1$.
%If $\nabla_{\mathbf v} p_0 \in  L^\infty_{\mathbf x}({\mathbb R}^N,
%L^1_{\mathbf v}( {\mathbb R}^N)) $, then
%$\nabla_{\mathbf v} p \in  L^{\infty}(0,T;L^\infty_{\mathbf x}(\Omega,
%L^1_{\mathbf v}( {\mathbb R}^N)) $ and the solution is unique.
\end{theo}

{\bf Proof.}  The solution is constructed as the limit of solutions
$p_{m}$  of approximating problems defined
by equations (\ref{lin:ec})-(\ref{lin:in}) with boundary
condition of the form:
\begin{eqnarray}
p_{m}^-({\mathbf x},{\mathbf v},t) = 
g(p_{m-1}^+({\mathbf x},{\mathbf v},t))
\quad {\rm on} \quad \Sigma_T^-,
\label{bc:rpmeps}  
\end{eqnarray}
The operators $g$ defining these boundary conditions in formulas
(\ref{bc:Nr0})-(\ref{bc:Nr1}) are positive.
The proof follows the same lines as the proof of 
Theorem \ref{maintheorem}, with changes to handle the boundary 
conditions, that we summarize.

The scheme is well defined  thanks to Theorems \ref{positivity} 
and \ref{conservation}, starting from $p_1=0$ and choosing
boundary values $p_2^-$ for $p_2$ on $\Sigma_T^-$
with the regularity (\ref{bpm}).
All the estimates established in Step 3 of the proof of Theorem 
\ref{maintheorem} hold. However, now $g=g(p_{m-1}^+({\mathbf x},
{\mathbf v},t))$
and we need to obtain uniform bounds of the traces at the boundaries.
Let us analyze the explicit expressions given by (\ref{bc:Nr0})-(\ref{bc:Nr1}).

We analyze first the bounds associated to the boundary condition at $S_{r_0}.$
From identity (\ref{bc:Nr0}), we deduce:
\begin{eqnarray}
 p^-_{m}(r_0, \boldsymbol \theta,v_r, \boldsymbol \phi,t) \!=\!
 {e^{-{\beta \over \sigma}|\mathbf v - \mathbf v_0|^2} \over {\cal I}_0} \Big[ 
 \int_0^\infty  \hskip -4mm  d \tilde v_r  \tilde v_r^{N-1} \hskip -2mm
 \int_{ \{\tilde {\boldsymbol \phi} \in S_{N\!-\!1} 
 | \tilde{\mathbf v} \cdot  {\mathbf n} <0 \} } 
 \hskip -18mm d \tilde {\boldsymbol \phi} \, 
 p^-_{m-1}(r_0,\boldsymbol \theta,\tilde v_r, \tilde {\boldsymbol \phi},t)
 \Big]. \label{bc0-e1}
\end{eqnarray} 
Multiplying (\ref{bc0-e1}) by $v_r^{N-1}$ and integrating over 
$\Sigma_T^-$, we find:
\begin{eqnarray}
 \|p^-_{m} \|_{L^1(\Sigma^-_T)} =
 \|p^-_{m-1} \|_{L^1(\Sigma^-_T)}
 = \ldots =  \|p^-_{2} \|_{L^1(\Sigma^-_T)}.
 \label{bc0-e2}
\end{eqnarray} 
Multiplying (\ref{bc0-e1}) by $v_r^{N-1} v_r^\ell$, $\ell >0$, 
integrating,  and inserting (\ref{bc0-e2}) we obtain:
\begin{eqnarray}
\| |\mathbf v|^\ell p^-_{m} \|_{L^1 (\Sigma^-_T)} \leq 
\int_{0}^\infty \hskip -2mm d v_r \; v_r^{N-1+\ell}  
\int_{ \{ \tilde{ \boldsymbol  \phi} \in S_{N\!-\!1} 
| \tilde{\mathbf v} \cdot   {\mathbf n} <0 \} } 
 \hskip -12mm d \tilde {\boldsymbol \phi} \;
 {e^{-{\beta \over \sigma}
  |\mathbf v - \mathbf v_0|^2} \over {\cal I}_0 } \;
 \|p^-_{2}\|_{L^1(\Sigma^-_T)}. \label{bc0-e3}
\end{eqnarray} 
Multiplying (\ref{bc0-e1}) by $v_r^{N-1}$ and integrating over 
$(0,\infty) \times 
\{  \tilde {\boldsymbol \phi} \in S_{N\!-\!1}
 |  \tilde {\mathbf v} \cdot   {\mathbf n} <0 \} $, 
we find:
\begin{eqnarray}
 \int_0^\infty  \hskip -4mm  d v_r  v_r^{N-1} \hskip -2mm
 \int_{ \{ {\boldsymbol \phi} \in S_{N\!-\!1} | {\mathbf v} \cdot {\mathbf n} <0 \} }
 \hskip -19mm d {\boldsymbol \phi} \,
 p^-_{m}(r_0, \boldsymbol \theta,v_r, \boldsymbol \phi,t) \!=\!
 \int_0^\infty  \hskip -4mm  d \tilde v_r  \tilde v_r^{N-1} \hskip -2mm
 \int_{ \{\tilde {\boldsymbol \phi} \in S_{N\!-\!1} 
 | \tilde{\mathbf v} \cdot  {\mathbf n} <0 \} } 
 \hskip -19mm d \tilde {\boldsymbol \phi} \, 
 p^-_{m-1}(r_0,{\boldsymbol \theta},\tilde v_r, 
 \tilde {\boldsymbol \phi},t) \nonumber \\
= \ldots = 
  \int_0^\infty  \hskip -4mm  d \tilde v_r  \tilde v_r^{N-1} \hskip -2mm
 \int_{ \{\tilde {\boldsymbol \phi}
  \in S_{N\!-\!1} | \tilde{\mathbf v} \cdot   {\mathbf n} <0 \} } 
 \hskip -19mm d \tilde {\boldsymbol \phi} 
 \, p^-_{2}(r_0, \boldsymbol \theta,\tilde v_r, \tilde 
 {\boldsymbol \phi},t).
\end{eqnarray} 
Therefore, for any $q \in (1, \infty)$:
\begin{eqnarray}
 \|p^-_{m} \|_{L^q(\Sigma^-_T)}^q  \!\leq \!
\int_{0}^\infty \hskip -4mm d v_r \hskip -1mm
\int_{ \{ {\boldsymbol \phi} \in S_{N\!-\!1} | {\mathbf v} \cdot   {\mathbf n} <0 \} } 
 \hskip -20mm d  {\boldsymbol \phi} \;  v_r^{N-1}  
 {e^{-{q \beta \over \sigma}
  |\mathbf v - \mathbf v_0|^2} \over {\cal I}_0^q }
 \| \hskip -1mm \left[ \hskip -1mm
 \int_0^\infty  \hskip -4mm  d \tilde v_r  \tilde v_r^{N\!-\!1} \hskip -2mm
 \int_{ \{\tilde \phi \in S_{N\!-\!1} | \tilde{\mathbf v} \cdot   {\mathbf n} <0 \} } 
 \hskip -21mm d \tilde \phi \, p^-_{2}(r_0,\theta,\tilde v_r, \tilde \phi,t)
 \hskip -1mm \right]^q \hskip -3mm \|_{L^\infty_{\boldsymbol \theta,t}}.
%\label{bc0-e5}
\nonumber
\end{eqnarray}
We may estimate uniformly the norms
$\|p^-_{m} \|_{L^\infty(\Sigma^-_T)}$, 
$\|(1+|\mathbf v|^2)^{\mu\over 2} p^-_{m} \|_{L^\infty(\Sigma^-_T)}$
and $\| p^-_{m} \|_{L^q_k(\Sigma^-_T)}$, for $1 \leq q \leq \infty,$  
in a similar way. 

Let us recall the boundary condition at $r=r_1$:
\begin{eqnarray}
 p^-_m(r_1,\! \boldsymbol \theta, \! v_r, \! \boldsymbol \phi, \!t) \!=\!  
 {e^{-{\beta \over \sigma} |\mathbf v \!-\! \mathbf v_0|^2} \over |{\cal I}_1|}
 \Big[  j_0 \!-\hskip -2mm  \int_0^\infty \hskip -4mm
d \tilde v_r   \tilde v_r^{N-1} \hskip -2mm
\int_{ \{ \tilde {\boldsymbol \phi} \in S_{N\!-\!1} | 
\tilde {\mathbf v} \cdot \mathbf n  >0 \} } 
\hskip -22mm
d \tilde {\boldsymbol \phi} \, p^+_{m-1}(r_1, \! \boldsymbol \theta, \! \tilde v_r, 
\!\tilde {\boldsymbol \phi}, \!t) f_1(\tilde{\mathbf v})  \Big]. 
\label{bc:Nr1eps}
\end{eqnarray}
Multiplying by $\mathbf v \cdot \mathbf n$ we get:
\begin{eqnarray}
\|p_m^-\|_{L^\infty_k(\Sigma_T^-)} \!\leq\! 
K_1(\beta,\sigma,\chi,\sigma_v)
\left[ \|j_0\|_\infty \!+\!  
K_2(\chi,\sigma_v) \|  p^+_{m\!-\!1} \|_{L^\infty_k(\Sigma_T^+)} 
\right].
\label{bc1-e0}
\end{eqnarray}

%$|\hat{\mathbf v}| e^{-|\hat{\mathbf v} |^2}$ that is below
%$0.45 <1$. Due to the factor $\hat{\mathbf v} \!\cdot\! \mathbf n$, the
%quotient of the integrals should be smaller than one too.

Set $\omega_q={N\beta \over q'} + \|a^-\|_\infty $.
From identity (\ref{lp}) in Theorem \ref{conservation}, we deduce:
\begin{eqnarray}
\|e^{-\omega_q t}
p_{m-1}^+\|_{L^q_k(\Sigma_T^+)} \leq 
\Big[ \| p_0Ê\|_{L^q(\Omega \times \mathbb R^N)} + 
\|e^{-\omega_q t} p_{m-1}^-\|_{L^q_k(\Sigma_T^-)}\Big],
\label{bc1-e1}
\end{eqnarray}
for any $q \in [1,\infty].$ Set $\omega = N\beta  + \|a^-\|_\infty.$
Multiplying equation (\ref{bc:Nr1eps}) by $e^{-\omega t}{\mathbf  v} \cdot  
\mathbf n$ and integrating over $\Sigma_T^-$, we find that
\begin{eqnarray}
\|e^{-\omega t} p_{m}^-\|_{L^\infty_k(\Sigma_T^-)} \leq K_1
\| e^{-\omega t}j_0\|_\infty + 
K_1K_2 \|e^{-\omega t} p_{m-1}^+\|_{L^\infty_k(\Sigma_T^+)}.
\label{bc1-e2}
\end{eqnarray}
Inserting (\ref{bc1-e1}) in (\ref{bc1-e2}) and iterating we obtain:
\begin{eqnarray}
\|e^{-\omega t} p_{m}^-\|_{L^\infty_k(\Sigma_T^-)} \!\leq\!  {1\over 1\!- \!K_1K_2}
\left[ C(j_0,\!\omega,\!T,\!K_1) 
\!+\! \| p_0Ê\|_{L^\infty_{\mathbf x \mathbf v}} \right]
\nonumber \\
\!+\! (K_1K_2)^{m-2} \|e^{-\omega t} p_{2}^-\|_{L^\infty_k(\Sigma_T^-)}.
\label{bc1-e3}
\end{eqnarray}
Using (\ref{bc1-e1}), we extend this uniform estimate to 
$\| p_{m}^+\|_{L^\infty_k(\Sigma_T^+)}.$

Multiplying equation (\ref{bc:Nr1eps}) by $|{\mathbf  v} |^{\ell}$, $\ell=0,\ldots,\mu$,
we find:
\begin{eqnarray}
\| |\mathbf v |^\ell p_m^-\|_{L^\infty(\Sigma_T^-)} \!\leq\!  { \| |\mathbf v |^\ell
e^{-{\beta \over \sigma} |\mathbf v \!-\! \mathbf v_0|^2}\|_{\infty} \over |{\cal I}_1|}
\left[ \|j_0\|_\infty \!+\!  K_2 \;  \|  p^+_{m\!-\!1} \|_{L^\infty_k(\Sigma_T^+)} 
\right].
\label{bc1-e4} \\
\| |\mathbf v |^\ell p_m^-\|_{L^1(\Sigma_T^-)} \!\leq\!  { \| |\mathbf v |^\ell
e^{-{\beta \over \sigma} |\mathbf v \!-\! \mathbf v_0|^2}\|_{1} \over |{\cal I}_1|}
\; {\rm meas}(\Omega)
\left[ \|j_0\|_\infty \!+\!  K_2 \;  \|  p^+_{m\!-\!1} \|_{L^\infty_k(\Sigma_T^+)} 
\right].  \label{bc1-e5}
\end{eqnarray}
In a similar way, we bound  uniformly
$\| |\mathbf v |^\ell p_m^-\|_{L^\infty_k(\Sigma_T^-)}$
and  $\| |\mathbf v |^\ell p_m^-\|_{L^1_k(\Sigma_T^-)}$
for $\ell=0,\ldots,\mu$.

The above uniform estimates on the boundary values yield the
uniform estimates on $p_{m}$ in Steps 3 and 4 of Theorem 
\ref{maintheorem}. 
We can extract converging subsequences as
in Step 5, with $\mathbf F(c_m)$ and 
$a=b(p_m)+\alpha \nu$ fixed, and pass
to the limit in the weak formulation as in Step 6, with obvious 
simplifications.
For the boundary term, an extracted subsequence 
${\rm Tr \,} p^\pm_{m'} = p_{m'}^\pm \rightharpoonup \pm$
in $L^q(\Sigma_T^\pm)$ and $L^q_k(\Sigma_T^\pm)$ weak 
for $1 \leq q < \infty$ and weak*
for $q=\infty$. This allows to pass to the limit in the boundary
term but we must justify that $g^-$ and $g^+$ satisfy the
equations defining the boundary conditions. Multiplying
(\ref{bc:Nr0})-(\ref{bc:Nr1}) by a test function
 $\psi \in C_c(\Sigma_T)$ and integrating, we find
 \begin{eqnarray*}
 \int_{\Sigma_T^- \cap \{Ê|\mathbf x|=r_0\}} \hskip -16mm
 p^-_{m'} \psi dS d\mathbf v dt &\!=\!&  
  \int_{\Sigma_T^- \cap \{Ê|\mathbf x|=r_0\}} \hskip -16mm 
 e^{-{\beta \over \sigma}|\mathbf v - \mathbf v_0|^2} {\cal I}_0^{-1}  \left[
 \int_0^\infty  \hskip -4mm  d  \tilde v_r  \tilde v_r^{N-1}  \hskip -1mm
\int_{ \{\tilde {\boldsymbol \phi} \in S_{N\!-\!1} | \tilde{\mathbf v} 
\cdot \mathbf n  <0 \} } 
 \hskip -16mm d \tilde   {\boldsymbol \phi} \; p^-_{m'-1}  \right] \, 
 \psi dS d\mathbf v dt, \; \label{bc:Nr0int} \\
 \int_{\Sigma_T^- \cap \{Ê|\mathbf x|=r_1\}}  \hskip -16mm
 p^-_{m'}   \psi dS d\mathbf v dt &\!=\!&  
 \int_{\Sigma_T^- \cap \{Ê|\mathbf x|=r_1\}} \hskip -16mm 
 e^{-{\beta \over \sigma} |\mathbf v \!-\! \mathbf v_0|^2} 
j_0   \, \psi dS d\mathbf v dt
 \nonumber \\
+&&  \hskip -8mm
 \int_{\Sigma_T^- \cap \{Ê|\mathbf x|=r_1\}}  \hskip -16mm
 e^{-{\beta \over \sigma} |\mathbf v \!-\! \mathbf v_0|^2}  |{\cal I}_1|^{-1}   
 \left[ \int_0^\infty \hskip -4mm
d \tilde v_r   \tilde v_r^{N-1} \hskip -1mm
\int_{ \{ \tilde {\boldsymbol \phi} \in S_{N\!-\!1} | 
\tilde {\mathbf v} \cdot \mathbf n  >0 \} } 
\hskip -16mm
d \tilde {\boldsymbol \phi} \, p^+_{m'-1}  f_1(\mathbf {\tilde v})
 \right] \, \psi dS d\mathbf v dt. \quad \label{bc:Nr1int}
\end{eqnarray*}      
Taking limits, the same identities  hold for $g^+$ and $g^-$.\\
           
Once we have understood the difficulties introduced by the
nonlocal boundary conditions, we can combine the strategies
developed in the proofs of Theorems \ref{maintheorem}
and \ref{maintheorembc} to obtain an existence result for
the original angiogenesis problem.

\begin{theo} 
\label{maintheoremfull}
Let us assume that
\begin{eqnarray}
p_0 \geq 0, c_0 \geq 0, \label{signofull} \\
c_0  \in W^{2,\infty}(\Omega),  \label{c0inffull} \\
(1+|{\mathbf v}|^2)^{\mu/2} p_0
\in L^{\infty} \cap L^1(\Omega \times {\mathbb R}^N), \quad \mu >N,
\label{p0infull} \\
c_{r_0}\in L^{\infty}(0,T; L^{\infty}(S_{r_0})), \label{bcmfull} 
%c_{r_0}\in L^{\infty}(0,T; W^{2,\infty}(S_{r_0})), \label{bcmfull} 
\end{eqnarray}
and that a function $c_b$ is found verifying the hypotheses of Theorem \ref{estimatesc}.
Then, there exists a positive solution $(p, c)$ of the initial value problem 
(\ref{eq:c})-(\ref{bc:c}) with  boundary conditions given by 
(\ref{bc:Nr0})-(\ref{bc:Nr1}) satisfying:
\begin{eqnarray}
c \in L^{\infty}(0,T;W^{1,\infty}(\Omega)), \label{dcinffull} \\
p \in L^{\infty}(0,T;L^{\infty}\cap L^1(\Omega \times {\mathbb R}^N)),
\nabla_{\mathbf v }p \in L^{2}(0,T;L^2(\Omega \times {\mathbb R}^N)),
\label{pinffull} \\
(1+|{\mathbf v}|^2)^{\mu/2}p \in 
L^{\infty}(0,T;L^{\infty}\cap L^1(\Omega \times {\mathbb R}^N)),  
\label{pdecayfull} \\
(1+ |\mathbf v \cdot \mathbf n|)(1+|{\mathbf v}|^2)^{\mu/2} {\rm Tr \,}p^\pm \in 
L^{\infty}(0,T;L^{\infty}\cap L^1(\Sigma_T^\pm)),  
\label{ptracefull} \\
p \in L^{\infty}(0,T;L^\infty_{\mathbf x}(\Omega,
L^1_{\mathbf v}( {\mathbb R}^N)), \label{pmominffull}  
\end{eqnarray}
provided the functions $K_1,K_2$ defined in (\ref{K1K2})
satisfy $K_1K_2<1$.
The norms of the solution are bounded in terms of the norms
of the data and the parameters.
\end{theo}

{\bf Proof.}

We consider the scheme (\ref{eq:pm})-(\ref{bc:cm}) with boundary
conditions (\ref{bc:rpmeps}), where $g$ is given by 
(\ref{bc:Nr0})-(\ref{bc:Nr1}).
We set $p_1=0$, so that $c_1$ is the solution of a heat equation.
Then, $p_2$ is the solution of the problem with bounded
coefficients ${\mathbf F(c_1)}$ and $\alpha(c_1)$ and fixed
boundary condition $p_2^-$ with the regularity (\ref{bpm}).
As in step 1 of the Proof of Theorem \ref{maintheorem}, the
sequence of iterates $(c_{m-1},p_m)$ is well defined thanks to
Proposition \ref{comparisonneuman}, Theorem \ref{estimatesc},
Theorems \ref{positivity}, \ref{conservation} and \ref{comparison},
Lemma \ref{cotaj} and Proposition \ref{velocityLinf}. 
The iterates are nonnegative, and the coefficients $j(p_{m-1})$,
$b(p_{m-1})$, $\alpha(c_{m-1})$ and ${\mathbf F}(c_{m-1})$
are bounded functions. As we have seen in the proof of Theorem
\ref{maintheorembc}, the boundary conditions for $p_m^-$
satisfy the regularity (\ref{bpm}).

The estimates for $c_m$ and $p_m$ in Steps 2 and 3 of the
Proof of Theorem \ref{maintheorem} hold. However, we do not
obtain immediate uniform estimates on the $L^q$ norms
of $p_m$ unless we estimate first the boundary conditions.
Setting $a=b(p_{m-1}) - \alpha(c_{m-1}) \nu$, we
have $\|a^-\|_\infty \leq \alpha_1 \| \nu \|_{\infty}$. Then, we may reproduce
the computations in the Proof of Theorem \ref{maintheorembc}
to get uniform bounds of $(1+ |\mathbf v \cdot \mathbf n|)
(1+|\mathbf v|^2)^{\mu/2} p_m^-$ in $L^1\cap L^\infty(\Sigma_T^+)$.
This provides uniform estimates on the $L^q$ norms of $p_m$
thanks to Theorem \ref{conservation}.
Steps 4, 5 and 6 proceed as in the proof of Theorem \ref{maintheorem}.
The passage to the limit in the boundary conditions is analogous to
that in the proof of Theorem \ref{maintheorembc}. The final solution
inherits all the bounds established for the iterates, as a result
of weak convergences.

\vskip 1cm

{\bf Acknowledgements.} This work has been supported by 
MINECO grant No. MTM2014-56948-C2-1-P.
The authors thank LL Bonilla and V. Capasso for suggesting the 
problem and for insight on the modeling.

\end{document}